\documentclass[reqno,english]{amsart}
\usepackage{amsfonts,amsmath,latexsym,verbatim,amscd,mathrsfs,color,array}
\usepackage[colorlinks=true]{hyperref}

\usepackage{amsmath,amssymb,amsthm,amsfonts,graphicx,color}
\usepackage{amssymb}
\usepackage{pdfsync}
\usepackage{epstopdf}
\usepackage[colorlinks=true]{hyperref}

\makeatletter
\def\@currentlabel{2.1}\label{e:dispaa}
\def\@currentlabel{2.21}\label{e:dispau}
\def\@currentlabel{2.22}\label{e:dispav}
\def\@currentlabel{2.23}\label{e:dispaw}
\def\@currentlabel{2.24}\label{e:dispax}
\def\theequation{\thesection.\@arabic\c@equation}
\makeatother

\hypersetup{linkcolor=black,citecolor=black,filecolor=black,urlcolor=black}

\newcommand{\inn}{{\quad\hbox{in } }}
\newcommand{\ttt}{\tilde }

\newcommand{\QQ}{{\mathcal Q}  }
\newcommand{\nn}{ {\nabla}  }

\newcommand{\by}{ {\bf y } }

\newcommand{\py}{ {\tt y } }
\newcommand{\bN}{ {\bf N } }
\newcommand{\hh}{ {\tt h } }
\newcommand{\vv}{ {\tt v } }

\newcommand{\pp}{ {\partial} }
\newcommand{\ww}{{\tt w}  }
\newcommand{\HH}{{{\mathbb H}}  }

\newcommand{\vp}{\varphi}

\newcommand{\RR}{{{\mathbb R}}}

\newcommand{\NN}{ {\mathcal N}}

\newcommand{\R} {\mathbb R}
\newcommand{\cuad}{{\sqcap\kern-.68em\sqcup}}

\newcommand{\BB}{{\tt B}}
\newcommand{\NNN}{{\tt N}}

\newcommand{\KK}{{\mathcal K}}

\newcommand{\foral}{\quad\mbox{for all}\quad}
\newcommand{\ve}{\varepsilon}

\newcommand{\be}{\begin{equation}}
\newcommand{\ee}{\end{equation}}

\newcommand{\la}{\lambda}

\newcommand{\equ}[1]{(\ref{#1})}

\renewcommand{\theequation}{\thesection.\arabic{equation}}
 
 \newtheorem{lemma}{Lemma}[section]

\newtheorem{theorem}{Theorem}
\newtheorem{prop}{Proposition}[section]
\newtheorem{corollary}{Corollary}[section]
\newtheorem{remark}{Remark}[section]
\newcommand{\bremark}{\begin{remark} \em}
\newcommand{\eremark}{\end{remark} }

\makeatother

\def\bb{\begin{equation}}
\def\bbs{\begin{equation*}}
\def\ees{\end{equation*}}

\def\ba{\begin{eqnarray}}
\def\ea{\end{eqnarray}}
\def\bas{\begin{eqnarray*}}
\def\eas{\end{eqnarray*}}

\newcommand{\Ga}{\Gamma}
\newcommand{\ga}{\gamma}

\newcommand{\bR}{{\bf R}}
\newcommand{\bS}{{\bf S}}

\newcommand{\MM}{{\mathcal M}}

\renewcommand{\theequation}{\arabic{section}.\arabic{equation}}

\begin{document}

\title[interface foliation on Riemannian manifolds]
{interface foliation near minimal submanifolds in Riemannian
manifolds with positive Ricci curvature}

\author[M. del Pino]{Manuel del Pino}
\address{\noindent M. del Pino - Departamento de
Ingenier\'{\i}a  Matem\'atica and CMM, Universidad de Chile,
Casilla 170 Correo 3, Santiago,
Chile.}
\email{delpino@dim.uchile.cl}

\author[M. Kowalczyk]{Michal Kowalczyk}
\address{\noindent  M. Kowalczyk - Departamento de
Ingenier\'{\i}a  Matem\'atica and CMM, Universidad de Chile,
Casilla 170 Correo 3, Santiago,
Chile.}
\email{kowalczy@dim.uchile.cl}

\author[J. Wei]{Juncheng Wei}
\address{\noindent J. Wei - Department of Mathematics, Chinese University of Hong Kong, Shatin, Hong Kong.}
\email{wei@math.cuhk.edu.hk}

\author[J. Yang]{Jun Yang}
\address{\noindent J. Yang - College of Mathematics and Computational Sciences,
  Shenzhen University,
  Nanhai Ave 3688, Shenzhen,  China, 518060. }
  \email{jyang@szu.edu.cn}

\begin{abstract}
Let $(\MM ,{\tilde g})$ be an $N$-dimensional smooth compact Riemannian manifold. We
consider the singularly perturbed Allen-Cahn equation
$$
\ve^2\Delta_{ {\tilde g}} {u}\,+\, (1 - {u}^2)u \,=\,0\quad \mbox{in } \MM,
$$
where $\ve$ is a small parameter.  Let $\KK\subset \MM$ be an
$(N-1)$-dimensional smooth minimal submanifold that separates $\MM$ into two disjoint components. Assume that
$\KK$ is non-degenerate in the sense that it does not support
non-trivial  Jacobi fields, and that $|A_{\KK}|^2+\mbox{Ric}_{\tilde g}(\nu_{\KK}, \nu_{\KK})$ is positive along $\KK$.
Then for each integer $m\geq 2$, we establish the existence of a sequence $\ve = \ve_j\to 0$,
and solutions $u_{\ve}$ with $m$-transition layers near $\KK$, with mutual distance $O(\ve|\ln \ve|)$.



\end{abstract}

\date{}\maketitle

\tableofcontents

\section{Introduction} \label{section1}
\setcounter{equation}{0}

In the gradient theory of phase transitions by Allen-Cahn
\cite{ac}, two phases of a material, $+ 1$ and $-1$  coexist in a
region $\Omega\subset \R^N$ separated by an $(N-1)$-dimensional
interface. The phase is idealized as a smooth $\ve$-regularization
of the discrete function, which is selected as a critical point of
the energy
$$
I_\ve(u) =  \int_\Omega \frac \ve 2 |\nabla u|^2   + \frac 1{4\ve} (1-u^2)^2,
$$
where $\ve>0$ is a small parameter. While any function with values
$\pm 1$ minimizes exactly the second term, the presence of the
gradient term conveys a balance in which the interface is selected
asymptotically as stationary for perimeter.  The energy $I_\ve$
may be regarded as an {\em $\ve$-relaxation} of the surface area:
indeed, in \cite{modica1} it is established that a sequence of
local minimizers $u_\ve$, with uniformly bounded energy, must
converge in $L^1_{loc}$-sense to a function of the form $\chi_E -
\chi_{E^c}$  so that $\partial E$ locally minimizes perimeter,
thus being a (generalized) minimal surface. This is the
starting point of the {\em $\Gamma$-convergence theory}, in which
the constraint of $I_\ve$ to a suitable class of separating-phase
functions, converges to the perimeter function of the interface.
Indeed, analogous assertions hold true for general families of
critical points, and for stronger notions of interface
convergence, see \cite{caffarellicordoba,pt,t}. The principle
above applies to modeling phase transition phenomena in many
contexts:  material science, superconductivity, population
dynamics and biological pattern formation, see for instance
\cite{js} and references therein.

\medskip
It is natural to consider situations in which phase transitions
take place in a manifold rather than in a subset of Euclidean
space.
In this paper we consider a compact  $N$-dimensional Riemannian
manifold $(\MM,{\tilde g})$, and want to investigate critical
points in $H^1(\MM)$ of the functional
$$
J_\ve(u) =  \int_\MM \frac{\ve}{2} \big|\nabla_{\tilde g}{\tilde u}\big|^2   + \frac 1{4\ve} \big(1-{\tilde u}^2\big)^2,
$$
with sharp transitions between $-1$ and $1$ taking place near
a $(N-1)$-dimensional minimal submanifolds of $\MM$. Critical points
of $J_\ve$ correspond precisely to classical solutions of the
Allen-Cahn equation in $\MM$,
\begin{equation}\label{originalproblem}
\ve^2\Delta_{ {\tilde g}} {\tilde u}\,+\, \big(1-{\tilde u}^2\big){\tilde u} \,=\,0\quad \mbox{in } \MM,
\end{equation}
where  $\Delta_{\tilde g}$ is the Laplace-Beltrami operator on
$\MM$.

\medskip
We let in what follows $\KK$  be a minimal $(N-1)$-dimensional
embedded submanifold of $\MM$,  which divides $\MM$ into two open
components $\MM_\pm$.  (The latter condition is not needed in some cases.) The {\em Jacobi operator} $\mathcal{J}$ of
$\KK$, corresponds to the second variation of $N$-volume along
normal perturbations of $\KK$ inside $\MM$: given any smooth small
function $v$ on $\KK$, let us consider the manifold $\KK(v)$, the
normal graph on $\KK$ of the function $v$, namely the image of
$\KK$ by the map $p\in\KK\mapsto
\mbox{exp}_p\bigl(v(p)\nu_{\KK}(p)\bigr)$. If $H(v)$ denotes the
mean curvature of $\KK(v)$, defined as the arithmetic mean of the
principal curvatures, then the linear operator $\mathcal{J}$ is
the differential of the map $v \mapsto nH(v)$ at $v=0$. More
explicitly, it can be shown that

\begin{equation}
\mathcal{J}\psi \ =\ \Delta_{\KK}\psi +\, |A_{\KK}|^2\psi
+\mbox{Ric}_{\tilde g}(\nu_{\KK},\nu_{\KK}) \psi , \label{jac}\ee
where $\Delta_{\KK}$ is the Laplace-Beltrami operator on
$\KK$, $|A_{\KK}|^2$ denotes the norm of the second fundamental
form of $\KK$, $\mbox{Ric}_{\tilde g}$ is the Ricci tensor of
$\MM$ and $\nu_{\KK}$ is a unit normal to $\KK$. We will briefly review these concepts
in Section 2.

\medskip
The minimal submanifold $\KK$ is said to be {\it nondegenerate} if
the are no nontrivial smooth solutions to the
 homogeneous problem
\begin{align}\label{Jacobiequation}
\mathcal{J}\psi=0\quad \hbox{in } \KK.
\end{align}
This condition implies that $\KK$ is isolated as a minimal
submanifold of $\MM$.

\medskip
In \cite{pacardritore}, Pacard and Ritor\'e assume that $\KK$ is
non-degenerate and, and proved that there exists a solution $u_\ve$ to
equation \equ{originalproblem} with values close to $\pm 1$ inside
$\MM_\pm$, whose (sharp) 0-level set is a smooth manifold which
lies $\ve$-close to $\KK$. More precisely, let  $w(z) := \tanh
\left (\frac z{\sqrt{2}}\right )$ be the unique solution of the
problem \bb w''+w- w^3=0\quad \mbox{in } \mathbb{R},\
w(0)=0,\ w(\pm\infty)=\pm 1,\label{definitionofH} \ee and
denote by $c_*$ its total energy, namely
$$
c_*:= \int_\R   \frac 12 |w'|^2 + \frac 14 (1- w^2)^2 .
$$
 Then the solution
$u_\ve$ in \cite{pacardritore} resembles near $\KK$ the function
$w(t/\ve)$, where $t$ is a choice of signed geodesic distance to
$\Gamma$. In particular
$$
J_\ve(u_\ve) \,\to\,  c_*\,|\KK|\,.\quad
$$

\medskip
In this paper we describe a new phenomenon induced by
the presence of positive curvature in the ambient manifold $\MM$:
in addition to non-degeneracy of $\KK$, let us assume that
\begin{align}\label{positivity}
K\,:=\,|A_{\KK}|^2+\mbox{Ric}_{\tilde
g}(\nu_{\KK},\nu_{\KK}) >0\quad\mbox{on }\KK.
\end{align}
Then, besides the solution by Pacard and Ritor\'e, there are
solutions with {\em multiple interfaces} collapsing onto $\KK$. In
fact, given any integer $m\ge 2$, we find a solution $u_\ve$ such that $u_\ve^2-1$
approaches $0$ in $\MM_\pm$ as $\ve \to 0$, with zero level
set constituted by $m$ smooth components $O(\ve|\log\ve|)$ distant
one to each other and to $\KK$, and such that
$$
J_\ve(u_\ve) \,\to\, m c_*\,|\KK|\,.\quad
$$
Condition (\ref{positivity}) is satisfied automatically if the
manifold $\MM$  has non-negative Ricci curvature. If $N=2$, $K$
corresponds simply to the Gauss curvature of $\MM$ measured along the
geodesic $\KK$.

\medskip
The nature of these solutions is drastically different from the
single-interface solution by Pacard and Ritor\'e\cite{pacardritore}.
They are
actually defined only if $\ve$ satisfies a {\em nonresonance
condition} in $\ve$. In fact, in the construction $\ve$ must
remain suitably away from certain values where a shift in Morse
index occurs. We expect that the solutions we find have a Morse
index $O(|\log\ve|^a)$ for some $a>0$
 as critical points of $J_\ve$, while the
single interface solution is likely to have its Morse index
uniformly bounded by the index of $\KK$ (namely the number of
negative eigenvalues of the operator $\mathcal{J}$).

\begin{theorem}\label{main}
Assume that $\KK$ is nondegenerate and embedded, and that condition
$(\ref{positivity})$ is satisfied. Then, for each $m\geq 2$, there
exists a sequence of values $\ve = \ve_j \to 0$
 such that problem $(\ref{originalproblem})$ has a
solution $u_{\ve}$  such that $u_\ve^2-1\to 0$  uniformly on compact
subsets of $\MM_{\pm}, $ while near $\KK$, we have
$$
u_\ve = \sum_{\ell =1}^m  w\left(\,\frac{{\tilde z} - \ve f_\ell ({\tilde y})}\ve\,
\right)\,+\,\frac{1}{2}\,\bigl((-1)^{m-1}\,-\,1\bigr) + o(1),
$$
where $({\tilde y},{\tilde z})$ are the Fermi coordinates defined near $\KK$
through the exponential map (see Section 2.1), and the functions $f_\ell$ satisfy
\begin{equation}
f_{\ell }({\tilde y})=  \Big(\ell -\frac{m+1}{2}\Bigr) \,
\, \left [ \, \sqrt{2}\log\frac{1}{{\ve}}
-\frac{1}{\sqrt{2}}\log\log\frac{1}{{\ve}}\right ] + O(1 ).
\end{equation}

\medskip
Moreover,
when $N=2$, there exist positive numbers $\nu_1,\ldots,\nu_{m-1}$ such that given $c>0$ and all
 sufficiently small $\ve>0$ satisfying
\begin{equation} \label{gap1}
\Big|\frac 1{\log\frac 1\ve}- \frac{\nu_i}{j^2}\Big|> \frac c {j^3},
\quad \foral\,  i=1,\ldots ,m-1,\quad j=1,2,\dots .
\end{equation}
a solution $u_\ve$ with the above properties exists.

\end{theorem}

\medskip
We observe that the same result holds if $m$ is {\em even} and $\MM\setminus\KK$
consists of just one component.  Thus the condition that $\KK$ divides $\MM$ into
two connected components is not essential in general.

\medskip
As we will see in the course of the proof, the equilibrium location of the interfaces is asymptotically governed by a small perturbation of
the  {\em Jacobi-Toda system}
\be
\ve^2\Big(\Delta_{\KK}f_j+\bigl(|A_{\KK}|^2+\mbox{Ric}_{\tilde g}(\nu_{\KK},\nu_{\KK})\bigr)f_j\Big)
\,-\, a_0 \,\bigl[\, e^{-(f_j -f_{j-1})}\,-\, e^{-(f_{j+1} - f_j)}\,\bigr]
\,=\,0, \label{toda}
\ee
 on $ \KK$, $j=1,\ldots, m$, with the conventions $f_0= -\infty,\ f_{m+1}=+\infty$.
 Heuristically, the interface foliation near ${\mathcal K}$ is possible due to  a balance
 between the interfacial energy, which  decreases as the interfaces approach
 each other, and the fact that the
 length or area of each individual interface increases as the interface is closer to $\KK$ since $\MM$
 is positively curved near $\KK$.

\medskip
 What is  unexpected, is the need of a {\em nonresonance condition } in order to solve the Jacobi-Toda system.
 A question which is of independent interest is the solvability of the Jacobi-Toda system without
 the condition (\ref{positivity}).
 Similar resonance has been observed  the problem
of  building foliations of a neighborhood of a geodesic  by CMC tubes considered in \cite{maz-pac0,maz_pac}. This has also been the case for (simple) concentration phenomena for various elliptic problems,  see  \cite{dkw,mahm,mm1,mm2}.

\medskip
 Our result deals with  situations in which the minimal submanifold is  locally but not globally  area
minimizing. In fact, since condition (\ref{positivity}) holds, the Jacobi operator has at least one negative eigenvalue, and
near $\KK$,  $\MM$ cannot have parabolic points.      In the case of a bounded domain $\Omega$ of
$\R^2$ under Neumann boundary conditions, a multiple-layer solution near a non-minimizing
straight segment orthogonal to the boundary was built in
\cite{dkw2}. In ODE cases  for the Allen-Cahn equation, clustering interfaces had been previously observed in
\cite{dy,na,nt}.
No resonance phenomenon is present in those situations, constituting a major qualitative difference with the current setting.

The method consists of linearizing the equation around the approximation
$$
u_0(x,z) =  \sum_{\ell =1}^m  w\left(\,\frac{{\tilde z} - \ve f_\ell ({\tilde y})}\ve\,
\right)\,+\,\frac{1}{2}\,\bigl((-1)^{m-1}\,-\,1\bigr),
$$
and then consider a projected form of the equation which can be solved boundedly
after finding a satisfactory linear theory, and then
applying the contraction mapping principe.
In that process the functions $f_j$ are left as arbitrary functions under some growth constraints.
At the last step one gets an equation which can be described as a small perturbation of the Jacobi-Toda system

\medskip
We do not expect that interface foliation
occurs  if the limiting interface is  a
minimizer of the perimeter since in such a case both perimeter of the interfaces
and their interactions decrease the energy, so no balance for their equilibrium locations
is  possible. On the other hand, negative Gauss curvature seems also prevent interface foliation.
This is suggested by a version of De Giorgi-Gibbons conjecture  for problem
\equ{originalproblem} with $\MM$ the hyperbolic space, established in \cite{birindelli 2}.

\section {Geometric background and the ansatz}\label{section2}
\setcounter{equation}{0}

{ In the first preliminary part of this section,  we list some necessary notions from differential geometry:
 Fermi coordinates near a submanifold
of $\MM$, minimal submanifold, as well as Laplace-Beltrami and Jacobi operators.
We then express the problem in a suitable form, define an approximate solution and estimate its error.

\subsection{Local coordinates }\label{section2.1}
Let $\MM$ be an $N\geq 2$-dimensional smooth compact Riemannian manifold without boundary with given  metric ${\tilde g}$.
We assume that $\KK$ is an $N-1$ dimensional submanifold of $\MM$.
For each given point $p\in\KK$,  $T_p\MM$ splits naturally as
$$
T_p\MM=T_p\KK\oplus N_p\KK,
$$
where $T_p\KK$ is the tangent space to $\KK$ and $N_p\KK$ is its normal complement, which spanned respectively by
orthonormal bases $\{E_i: i=1,\cdots,N-1\}$ and $\{E_N\}$. More generally, we have for the tangent and  normal bundles over $\KK$ the decomposition
$$
T\MM=T\KK\oplus N\KK.
$$
Let us  denote by $\bigtriangledown$ the connection induced by the metric ${\tilde g}$ and by $\bigtriangledown^{N}$
the corresponding normal connection on the normal bundle.

{\bf Notation:} Up to section \ref{2point4}, we shall always use the following convention for the  indices
$$i,j,k,l\cdots\in\{1,2,\cdots,N-1\}, \quad a,b,c,\cdots\in\{1,2,\cdots,N\}.$$

\medskip
Given $p\in\KK$, we use some geodesic coordinates ${\tilde y}$
centered at $p$. More precisely, in a neighborhood of $p$ in $\KK$, we consider normal geodesic coordinates
\begin{align}\label{Yp}
\ttt y= Y_p(\ttt \py)  :=\mbox{exp}_p^{\KK}({\tilde \py}_iE_i),\quad {\tilde \py}=({\tilde \py}_1,\cdots,{\tilde \py}_{N-1})\in {\mathcal V},
\end{align}
where $\mbox{exp}^{\KK}$ is the exponential map on $\KK$ and summation over repeated indices is understood. ${\mathcal V}$ is a neighborhood of the
origin in $\RR^{N-1} $.

\medskip
This yields the coordinate vector fields $X_i=f_*(\partial_i),\,i=1,\cdots,N-1$
where $f({\tilde \py})=Y_p(\ttt \py)$. For any $E\in T_p\KK$, the curve
$$s\rightarrow \ga_E(s)=\mbox{exp}_p^{\KK}(sE)$$
is a geodesic in $\KK$, so that
$$
\bigtriangledown_{X_i}X_j|_p\in N_p\KK\quad \mbox{for any }i,j=1,\cdots,N-1.
$$
We recall that the Christoffel symbols    $\Ga^N_{ij},\, i,j=1,\cdots,N-1$ are given by
$$
\nabla_{X_i} X_j|_p\,=\,\Ga^N_{ij}E_N,\quad\mbox{i. e.}\quad\Ga^N_{ij}\,=\,{\tilde g}(\nabla_{X_i}X_j,E_N)
.
$$
We also assume that at $p$ the normal vector $E_N$ is transported parallelly (with respect to $\bigtriangledown^N$)
through geodesics $\ga_E(s)$ from $p$. This yields a frame field $X_N$ for $N\KK$ in a neighborhood of $p$ which satisfies
$$
\nabla_{X_i} X_N|_p\in T_p\KK,\quad\mbox{i.e.}\quad{\tilde g}(\nabla_{E_i}E_N,E_N)|_p=0,\quad i=1,\cdots,N-1.
$$
We define the numbers $\Ga^j_{iN},\, i,j=1,\cdots,N-1,$ by
$$
\nabla_{X_i} X_N|_p\,=\,\sum_{j=1}^{N-1}\Ga^j_{iN}E_j,\quad\mbox{i.e.}\quad\Ga^j_{iN}\,=\,{\tilde g}(\nabla_{X_i} X_N,E_j).
$$

In a neighborhood of $p$ in $\MM$, we choose the {\it Fermi coordinates} $({\tilde \py}, {\tilde z})$ on $\MM$ defined by
\begin{align}
\Phi^0({\tilde \py},{\tilde z})=\mbox{exp}_{Y_p({\tilde \py})}({\tilde z}E_N)\quad \mbox{with }
({\tilde \py}, {\tilde z})=({\tilde \py}_1,\cdots,{\tilde \py}_{N-1},{\tilde z})
\in {\mathcal V} \times \bigl(-\delta_0,\delta_0\bigr),
\end{align}
where $\mbox{exp}_{ Y_p({\tilde \py}) }$ is the exponential map at $ Y_p({\tilde \py})$ in $\MM$.
We also have corresponding coordinate vector fields
$$
X_i=\Phi^0_*(\partial_{{\tilde y}_i}),\quad X_N=\Phi^0_*(\partial_{\tilde z}).
$$
By construction, $X_N|_p=E_N$.

\subsection{Taylor expansion of the metric}
In this section we will follow the notation and calculations of \cite{maz-pac0}.
By our choice of coordinates and the  Gauss Lemma, on $\KK$ the metric ${\tilde g}$ splits in the following way,
\begin{align}\label{gausslemma}
{\tilde g}(p)\,=\,\sum_{i,j=1}^{N-1}{\tilde g}_{ij}\mathrm{d}{\tilde y}_i\otimes\mathrm{d}{\tilde y}_j
\,+\,{\tilde g}_{NN}\mathrm{d}{\tilde z}\otimes\mathrm{d}{\tilde z},\quad p\in\KK.
\end{align}
As usual, the Fermi coordinates above are chosen so that the metric coefficients satisfy
$$
{\tilde g}_{ab}={\tilde g}(X_a,\,X_b)=\delta_{ab}\quad\mbox{at }\  p.
$$
Furthermore, ${\tilde g}(X_i,\,X_N)=0$ in some neighborhood of $p$ in $\KK$.
Then
$$
X_i{\tilde g}(X_j,\,X_N)={\tilde g}(\bigtriangledown_{X_i}X_j,\,X_N)+{\tilde g}(X_j,\,\bigtriangledown_{X_i}X_N)\quad\mbox{on }\KK,
$$
yield the identity
\begin{align}\label{identity}
\Ga^N_{ij}+\Ga^j_{iN}=0\quad\mbox{at }\ p.
\end{align}

We denote by $\Ga_i^j:N\KK\rightarrow \mathbb{R},\,i,j=1,\cdots,N-1,$ the 1-forms defined on the normal bundle of $\KK$ as
\begin{align}\label{christoffel}
\Ga_i^j(E_N)\,=\,{\tilde g}(\bigtriangledown_{E_i}E_j,E_N).
\end{align}
The  {\it second fundamental form} $A_{\KK}: T{\KK}\times T{\KK}\rightarrow N\KK$
of the submanifold $\KK$ and its corresponding norm are then given by
\begin{align}\label{secondfundamentalform}
A_{\KK}(E_i,E_j)=\Ga^j_i(E_N)E_N,\quad |A_{\KK}|^2=\sum^{N-1}_{i,j=1}\Bigl(\Ga^j_i(E_N)\Bigr)^2.
\end{align}
For $X,Y,Z,W\in T\MM$, the curvature operator and curvature tensor are respectively defined by the relations
\begin{align}\label{curturetensor}
&R(X,Y,Z)=\bigtriangledown_X\bigtriangledown_YZ-\bigtriangledown_Y\bigtriangledown_XZ-\bigtriangledown_{[X,Y]}Z,
\\
&R(X,Y,Z,W)={\tilde g}\bigl(R(Z,W)Y,X\bigr).
\end{align}
The {\it Ricci tensor} of $(\MM, {\tilde g})$ is defined by
\begin{align}\label{Rici}
\mbox{Ric}_{\tilde g}(X,Y)={\tilde g}^{ab}R(X,X_a,Y,X_b).
\end{align}

We now compute higher order terms in the Taylor expansions of the metric coefficients. The metric coefficients
at $q=\Phi^0(0,{\tilde z})$ are given in terms of geometric data at $p=\Phi^0(0,0)$
and $|{\tilde z}|=\mbox{dist}_{\tilde g}(p,q)$,
which is expressed by the next lemmas, see  Proposition 2.1 in \cite{maz-pac0} and the references therein.
\begin{lemma}\label{covariant}
At the point $q=\Phi^0(0,{\tilde z})$, the following expansions hold
\begin{align}
&\bigtriangledown_{X_N}X_N=O(|{\tilde z}|)X_a,
\\
&\bigtriangledown_{X_i}X_j=\Ga^j_i(E_N)X_N+O(|{\tilde z}|)X_a,\quad i,j=1,\cdots,N-1,
\\
&\bigtriangledown_{X_i}X_N=\bigtriangledown_{X_N}X_i=\sum_{j=1}^{N-1}\Ga^j_i(E_N)X_j+O(|{\tilde z}|)X_a,\quad i=1,\cdots,N-1.
\end{align}
\end{lemma}

\begin{lemma}
In the above coordinates $({\tilde y},{\tilde z})$, for any $i,j=1,2,\cdots,N-1$, we have
\begin{align}
{\tilde g}_{ij}(0,{\tilde z})&\,=\,\delta_{ij}-2\Ga^j_i(E_N){\tilde z}
-R(X_N,X_j,X_N,X_i)|{\tilde z}|^2\nonumber
\\
&\quad+\sum^{N-1}_{k=1}\Ga^k_i(E_N)\Ga^j_k(E_N)|{\tilde z}|^2+O(|z|^3),
\\
{\tilde g}_{iN}(0, z)&\,=\,O(|{\tilde z}|^2),
\\
{\tilde g}_{NN}(0,z)&\,=\,1+O(|{\tilde z}|^3).  
\end{align}
\end{lemma}

}
{
\subsection{The  Laplace-Beltrami and Jacobi operators}\label{section2.3}

If $(\MM,{\tilde g})$ is an $N$-dimensional Riemannian manifold, the
{\it Laplace-Beltrami operator} on $\MM$ is
defined in local coordinates by the formula
\begin{align}\label{laplace}
\Delta_{\MM}=\frac{1}{\sqrt{\mbox{det}\tilde g}}\,\partial_a\Bigl(\sqrt{\mbox{det}\tilde g}\, {\tilde g}^{ab}\partial_b\Bigr),
\end{align}
where ${\tilde g}^{ab}$ denotes the inverse of the matrix $({\tilde g}_{ab})$.
Let $\KK\subset \MM$ be an $(N-1)$-dimensional closed smooth embedded submanifold associated with the metric ${\tilde g}_0$ induced from $(\MM,{\tilde g})$.
Let $\Delta_{\KK }$ be the Laplace-Beltrami operator defined on $\KK$.

\medskip
Let us consider the space $C^{\infty}(N\KK)$,
which identifies with that of all smooth normal vector fields on $\KK$.
Since $\KK$ is a submanifold of codimension $1$, then  given a choice of orientation
and  unit normal vector field along $\KK$, denoted  by $\nu_{\KK}\in N\KK$,
we can write  $\Psi\in  C^{\infty}(N\KK)$ as $\Psi=\phi\nu_\KK$, where $\phi\in   C^{\infty}(\KK)$.

\medskip
For $\Psi\in C^{\infty}(N\KK)$, we consider the one-parameter family
of submanifolds $t\rightarrow \KK_{t,\Psi}$ given by
\begin{align}
\KK_{t,\Psi}\equiv\Bigl\{\mathrm{exp}_{\tilde y}\bigl(t\Psi({\tilde y})\bigr): {\tilde y}\in\KK\Bigr\}.
\end{align}
The first variation formula of the volume functional is defined as
\begin{align}\label{firstvariation}
\frac{\mathrm{d}}{\mathrm{d}t}\Big|_{t=0}\,\mbox{Vol}{(\KK_{t,\Psi})}=\int_{\KK}<\Psi,\,{\bf h}>_N\mathrm{d}V_{\KK},
\end{align}
where ${\bf h}$ is the mean curvature vector of $\KK$ in $\MM$, $<\cdot,\cdot>_N$ denotes the restriction
of $\tilde g$ to $N\KK$, and $\mathrm{d}V_{\KK}$ the volume element of $\KK$.

\medskip
The submanifold $\KK$ is said to be {\it minimal} if it is stationary point for the volume functional, namely if
\begin{align}
\frac{\mathrm{d}}{\mathrm{d}t}\Big|_{t=0}\,\mbox{Vol}{(\KK_{t,\Psi})}=0
\quad
\mbox{for any } \Psi\in C^{\infty}(N\KK),
\end{align}
or equivalently by (\ref{firstvariation}), if the mean curvature ${\bf h}$ is identically zero on $\KK$.
It is a standard fact that if $\Ga^j_i(E_N)$ is as in (\ref{christoffel}), then
\begin{align}\label{minimal}
\KK \mbox{ is minimal }\, \Longleftrightarrow\,\sum^{N-1}_{i=1}\Ga^i_i(E_N)=0.
\end{align}

\medskip
The {\it Jacobi operator} $\mathcal{J}$ appears in the expression of the second variation of the volume functional
for a minimal submanifold $\KK$
\begin{align}
\frac{\mathrm{d}^2}{\mathrm{d}t^2}\Big|_{t=0}\,\mbox{Vol}{(\KK_{t,\Psi})}
=-\int_{\Ga}<\mathcal{J}\Psi,\,\Psi>_N\mathrm{d}V_{\KK}
\quad
\mbox{for any } \Psi\in C^{\infty}(N\KK),
\end{align}
and is given by
\begin{align}\label{Jacobioperator}
\mathcal{J}\phi=-\Delta_{\KK}\phi-\mbox{Ric}_{\tilde g}(\nu_\KK, \nu_\KK)\phi-|A_{\KK}|^2\phi,
\end{align}
where $\Psi=\phi\nu_\KK$, as has been explained above.

\medskip
The submanifold $\KK$ is said to be {\it non-degenerate} if the Jacobi operator $\mathcal{J}$ is invertible,
or equivalently if the equation $\mathcal{J}\phi=0$ has only the trivial solution in $C^\infty(\KK)$

\color{black}
\subsection{Laplace-Beltrami Operator in Stretched Fermi Coordinates}\label{2point4}
To construct the approximation to a solution of (\ref{originalproblem}), which concentrates near $\KK$, after rescaling,
in $\MM/\ve$,
 we  introduce stretched Fermi coordinates in the neighborhood of the point $ \ve^{-1} p\in \ve^{-1} \KK$ by
\bb\label{definitionofzandz}
\Phi_{\ve}(\py,z)=\frac{1}{\ve}\Phi^0(\ve \py ,\ve z),\quad
(\py,z)=(\py_1,\cdots,\py_{N-1},z)\in  \ve^{-1}{\mathcal V}
\times  \Big(-\frac{\delta_0}{\ve},\frac{\delta_0}{\ve}\Big).
\ee
Obviously, in $\MM_{\ve}= \ve^{-1}\MM$ the new coefficients $g_{ab}$'s of the Riemannian metric,
after rescaling, can be written as
$$
g_{ab}(\py,z)={\tilde g}_{ab}(\ve \py, \ve z), \quad a,b=1,2,\cdots,N.
$$
\begin{lemma}
In the above coordinates $(\py,z)$, for any $i,j=1,2,\cdots,N-1$, we have
\begin{align}
{g}_{ij}(\py,z)&\,=\,\delta_i^j-2\ve\Ga^j_i(E_N)z-\ve^2R(X_j,X_N,X_N,X_i)|z|^2
\nonumber
\\
&\qquad+\ve^2\sum^{N-1}_{k=1}\Ga^k_i(E_N)\Ga^j_k(E_N)|z|^2+O(|\ve z|^3),
\\
{g}_{iN}(\py,z)&\,=\,O(|\ve z|^2),
\\
{g}_{NN}(\py,z)&\,=\,1+O(|\ve z|^3).
\end{align}
Here $R(\cdot)$ and $\Ga^{b}_a$ are computed at the point $p\in\KK$ parameterized by $(0,0)$.
\end{lemma}

Now we will  focus on the expansion of the Laplace-Beltrami operator defined by
\begin{align}
\begin{aligned}
\Delta_{\MM_\ve }&=\frac{1}{\,\sqrt{{\mbox{det}g}}\,}\,\partial_a
\Bigl(\,g^{ab}\,\sqrt{{\mbox{det}g}}\,\partial_b\,\Bigr)
\\
&=g^{ab}\,\partial_a\partial_b\,+\,(\partial_a g^{ab})\,\partial_b\,+\,\frac{1}{2}\,
\partial_a\bigl(\log{{(\mbox{det}g)}}\,\bigr)\,g^{ab}\,\partial_b.
\end{aligned}
\label{definitionofLaplacian0}
\end{align}
Using the assumption that the submanifold $\KK$ is minimal as in formula (\ref{minimal}), direct computation gives that
$$
\mbox{det} g=1-\ve^2K(\ve y)z^2+O(\ve^3|z|^3),
$$
where we have, using (\ref{secondfundamentalform}) and (\ref{Rici}), denoted
\begin{align}\label{definitionofK}
K=\mbox{Ric}_{\tilde g}(\nu_{\KK},\nu_{\KK})+|A_{\KK}|^2.
\end{align}
This gives
$$
\log{(\mbox{det} g)}=-\ve^2K(\ve y)z^2+O(\ve^3|z|^3).
$$
Hence, we have the expansion
\begin{align}
\begin{aligned}
\Delta_{\MM_\ve}=\pp_{zz} \,+\,\Delta_{\KK_\ve}                
\,+\,\ve^2 zK(\ve y)\,\pp_z  \,+\, B
\end{aligned}
\label{definitionofLaplacian}
\end{align}
where the operator $B$ has the form
\be
B\, =\, \ve z\,a^1_{ij}\,\pp_{ij} \,+\, \ve^2 z^2\, a_{iN}^2\, \pp_{iz} \,+\,  \ve^3z^3\,a^3_{NN}\,\pp_{zz}   \,+\, \ve^2 z\,b_i^1\pp_i   \,  +\, \ve^3 z^2\,b_N^2\,\pp_z\, .
\label{defB}\ee
and all the coefficients are smooth functions defined on a neighborhood of $\KK$ in $\MM$, evaluated at $(\ve y, \ve z)$.

}
\subsection{The local approximate solution} \label{section approz}
If we  set $u(x) := \tilde{u}(\ve x)$,  then problem (\ref{originalproblem}) is thus equivalent to
\begin{eqnarray}\label{scalingproblem}
 \Delta_{\MM_\ve }\, u\,+\,F(u) =0 \quad \hbox{ in }\MM_{\ve},
\end{eqnarray}
where $F(u) \,\equiv\,u\,-\,u^3$. In the sequel, we denote by $\MM_{\ve}$ and $\KK_{\ve}$ the $\ve^{-1}$-dilations of
 $\MM$ and $\KK$.

\medskip
To define the approximate solution we observe
 the heteroclinic solution to (\ref{definitionofH}) has the asymptotic properties
\begin{align}
\begin{aligned}
\label{propertiesofH}
w(z)\,-\,1&=\,- 2\,e^{-\sqrt{2}\,|z|}\,+\,O\bigl(e^{-2\sqrt{2}\,|z|}\bigr),\quad z>1,
\\
w(z)\,+\,1&=\, 2\,e^{-\sqrt{2}\,|z|}\,+\,O\bigl(e^{-2\sqrt{2}\,|z|}\bigr),\quad z<-1,
\\
w'(z)&=\, 2\,\sqrt{2}\, e^{-\sqrt{2}\,|z|}\,+\,O\bigl(e^{-2\sqrt{2}\,|z|}\bigr),\quad |z|>1,
\end{aligned}
\end{align}
For a fixed integer $m\geq 2$, we assume that the location of
the $m$ phase transition layers are characterized  in the coordinate $(y,z)$ defined in (\ref{definitionofzandz})  by the functions
$z\,=\,f_j(\ve y),\, 1\leq j\leq m$ with
$$ f_1(\ve y) < f_2(\ve y) < \cdots < f_m(\ve y),  $$  separated one to each other
by large distances  $O(|\log \ve|)$,
and define in coordinates $(y,z)$ the approximation
\be
u_0 (y,z) :=   \sum_{j=1}^m w_j \big(z- f_j(\ve y)\big)
\, +\,  \frac { (-1)^{m-1}-1} 2,
\quad
w_j(t):= (-1)^{j-1}w (t),
\label{ww}
\ee
with this definition we have that
$u_0(y,z) \approx w_j (z- f_j(\ve y))$ for values of $z$ close to $f_j(\ve y)$.

\medskip
The functions $f_j : \KK \to \R$ will be left as parameters,
on which we will assume a set of  constraints that we describe next.

\medskip
Let us  fix numbers $p> N$, $M>0$, and consider functions $h_j \in W^{2,p}(\KK)$, $j=1,\ldots, m$, such that

\be
\|h_j\|_{W^{2,p}(\KK)} := \|D^2_\KK h_j\|_{L^p(\KK)} + \|D_\KK h_j\|_{L^p(\KK)} + \|h_j\|_{L^\infty (\KK)}
\ \le \ M\, .
\label{assh}
\ee

For a small $\ve>0$, we consider the unique number $\rho= \rho_\ve $ with
\begin{equation}
e^{\,-\sqrt{2}\,\rho\,} =\, {\ve^2}\rho .
\label{rho}\end{equation}
We observe that $\rho_\ve $ is a large number that can be expanded in $\ve$ as
\begin{align*}
\rho_{\ve}\,=\,
\sqrt{2}\log\frac{1}{{\ve}}
-\frac{1}{\sqrt{2}}\log\Bigl(\sqrt{2}\log\frac{1}{{\ve}}\Bigr)
+O\Big(\frac{\log\log\frac{1}{{\ve}}}{\log\frac{1}{{\ve}}}\Big).
\end{align*}
Then we assume that the $m$ functions   $f_j : \KK \to \R$ are given by the relations
\be
f_k({\tilde y}) \ =\   \left (k\,- \, \frac {m + 1} 2 \right )\, \rho_\ve \,+ \, h_k({\tilde y}), \quad k= 1,\ldots, m,
\ee
so that
\be
 f_{k+1} ({\tilde y}) - f_{k}({\tilde y}) \,= \, \rho_\ve \,+ \,h_{k+1}({\tilde y}) - h_{k}({\tilde y})\, ,\quad k=1,2,\ldots, m-1.
\label{fj}
\ee
We will use in addition the conventions
$h_0\equiv -\infty$, $h_{m+1} \equiv +\infty$.

\medskip
Our first goal is to compute the error of approximation in a $\delta_0/\ve$ neighborhood of
$\KK_{\ve}$, namely the quantity:
\begin{eqnarray}
S(u_0) &\equiv&\Delta_{\MM_\ve }\,u_0\,+\,F(u_0).
\end{eqnarray}

\medskip
For each fixed $\ell$, $1\leq \ell \leq m$,
this error reproduces a  similar pattern on each set of the form
\begin{align}
A_\ell =\left\{\,(y,z)\in\KK_{\ve}\times \bigl(-\frac{\delta_0}{\ve},\frac{\delta_0}{\ve}\bigr)
\ \bigl / \   | z- f_\ell (\ve y) |\, \leq\, \frac 12 \rho_\ve + M \, \right \}.
\label{defan}
\end{align}
For $(y,z)\in A_\ell $,  we write $t= z- f_\ell (\ve y)$
and estimate in this range the quantity  $S(u_0)\big( y, t+ f_\ell (\ve y)\big)$.
We have the validity of the following expression.

\begin{lemma}\label{expansionerror}
For $\ell \in \{1,\ldots, m\}$ and $(y,z)\in A_\ell$ we have

\begin{align}
&(-1)^{\ell-1}S(u_0) (y,t+ f_\ell)\nonumber
\\
&\, =\,
 6(1- w^2(t))\ve^2 \rho_\ve \,\left  [ \, e^{- \sqrt{2}(h_\ell  -h_{\ell -1})} e^{ \sqrt{2} t}\,-\,
 e^{ -\sqrt{2}(h_{\ell +1}- h_\ell ) } e^{ -\sqrt{2} t}\, \right ] \label{expr2}
\\
&\quad\,-\, \ve^2\, \Bigl(\Delta_{\KK} h_\ell  +  (t + f_\ell  )K\,\Bigr) w'(t)
\,+\, \ve^2|\bigtriangledown_{\KK }h_\ell  |^2 w''(t) \, + \, (-1)^{\ell-1}\Theta_\ell(\ve y, t) \,.
\nonumber
\end{align}
where for some $\tau, \sigma>0$ we have
$$
\|\Theta_\ell (\cdot , t)\|_{L^p(\KK)} \le   C\ve^{2+ \tau}  e^{-\sigma|t|} .
$$

\end{lemma}

\proof

From (\ref{definitionofLaplacian}), using that $w_j'' + F(w_j)=0$,
we  derive that, for $(y,z)\in A_\ell$
\begin{align}
S(u_0)\, =\ &   F \big( u_0(y,z)\big)
\,-\,  \sum_{j=1}^m  F\big(w_{j} (z- f_j(\ve y)\big)
\,+\, \ve^2\sum_{j=1}^m|\bigtriangledown_{\KK }h_j (\ve y)|^2 w_{j}''\big(z- f_j(\ve y)\big)
\nonumber
\\
&\,-\, \ve^2\sum_{j=1}^m \,\Bigl(\Delta_{\KK} h_j (\ve y)+zK(\ve y)\,\Bigr) w_j'(z- f_j(\ve y))
\nonumber
\\
&\,+\,\ve^3z\, \big[\, a^1_{ik} (\ve y, \ve z) \,\pp_{ik} h_j(\ve y)   \,+ \,  b_i^1(\ve z, \ve y)\,\pp_i h_j (\ve y)\,\big]\, w_j'\big(z-f_j(\ve y)\big)
\nonumber
\\
\nonumber
\\
&\,+\,
\ve^3 \big[\, z^3\,a^3_{NN}(\ve y, \ve z)\,  +   z^2\, a_{iN}^2(\ve y, \ve z)\, \pp_{i} h_j (\ve y)\big] \,  w_j''\big(z-f_j (\ve y)\big)\,
\nonumber
\\
\nonumber
\\
& \,+\,
\ve^3 \,z\, a_{ik}^1(\ve y, \ve z) \,\pp_ih_j (\ve y)\,\pp_k h_j (\ve y)\,  w_j''\big(z-f_j (\ve y)\big)\,.
\label{expr1}
\end{align}
Let us consider first the case  $2\le \ell \le m-1$.

\medskip
We begin with the term
$$
F \big( u_0(y, t+ f_\ell)\big)
\, -\,  \sum_{j=1}^m  F\big(\,w_{j} ( t + f_\ell - f_j)\, \big),
\quad
|t|< \frac { \rho_\ve}2.
$$
 Since
$$
w(s) =   \pm \big(1 - 2 e^{-\sqrt{2}|s|} \big)
 \, +\,O \big( e^{-2\sqrt{2}|s|}\big)
 \quad
 \hbox{as } s\to \pm \infty ,
$$
we find that
for $j  < \ell$,
\be
w \big( t + f_\ell  - f_j\big) - 1  = -  2 e^{ -\sqrt{2}( f_\ell  - f_j)} e^{ -\sqrt{2} t} \,  +\,  O\left ( e^{ - 2\sqrt{2}|t + f_\ell  - f_j|}\, \right),
\label{r5}\ee
while for $j> \ell$,
\be
w \big( t + f_\ell  - f_j\big) +  1  =  2 e^{ -\sqrt{2}( f_j - f_\ell )} e^{ \sqrt{2} t} \,  +\,  O\left ( e^{ - 2\sqrt{2}|t + f_\ell  - f_j|}\, \right).
\label{r6}\ee

Now, since
\begin{align*}
F(w) = w(1-w^2),
\quad
|t| < \frac { \rho_\ve}2 + O(1),
\\
|f_\ell  - f_j| = |\ell-j|\rho_\ve +O(1),
\quad
e^{-\sqrt{2}\rho_\ve}= \ve^2 \rho_{\ve},
\end{align*}
we  find that if $|j-\ell|\ge 2$ and $0< \sigma < \sqrt{2}$, then for some $\tau>0$,

\be
|F( w_j( t + f_\ell  - f_j))|   \,\le \,C e^{ - \sqrt{2}|t + f_\ell  - f_j|}\, \le \,
 \ve^{2+\tau}\, e^{- \sigma |t|} .
\label{r1}\ee

On the other hand, for certain numbers $s_1,s_2\in (0,1)$ we have

\begin{align}
F\big( w ( t + f_\ell  - f_{\ell -1})\big)  =   F'( 1 ) \, a_1, + \frac 12 F''( 1+  s_1 a_1  ) \, a_1^2 ,
\label{r2}
\\
\nonumber
\\
F\big( w ( t + f_\ell  - f_{\ell +1})\big)  =   F'( 1)\, a_2 + \,\frac 12 F''( 1 -  s_2 a_2  ) \, a_2^2.
\label{r3}
\end{align}
where
$$
a_1 :=  w ( t + f_\ell  - f_{\ell -1}) -1,
\quad
a_2 :=  w ( t + f_\ell  - f_{\ell +1}) +1.
$$
Now,
we find
$$
(-1)^{\ell-1} u_0 =   w(t) - a_1 - a_2  -  a_3 , \quad a_3 = O\left ( \max_{|j- \ell|\ge 2} e^{ -\sqrt{2}|t + f_\ell  - f_j|}\, \right).
$$
Thus for some $s_3\in (0,1)$,
\begin{align}
(-1)^{\ell-1}F(u_0) \,=\,&  F(w)\,-\, F'( w)\, (a_1+a_2)
 \,+\,\frac 12 F''( w- s_3(a_1+a_2)) (a_1+a_2)^2\nonumber
\\
\nonumber
\\
& \,+\, O\left ( \max_{|j- \ell|\ge 2} e^{ -\sqrt{2}|t + f_\ell  - f_j|}\, \right).
\label{r4}
\end{align}
Combining relations \equ{r1}-\equ{r4} and using that
$$
F'(1)-F'(w) = 3(1-w^2),\quad   |a_1| + |a_2| =  O(e^{-\sqrt{2}\frac {\rho_\ve} 2})  =  O(e^{-\sqrt{2}|t|}),
$$
  we obtain
\begin{align*}
&(-1)^{\ell-1} \Big(\, F(u_0) - \sum_{j=1}^m F\big( w_j( t+ f_\ell  -f_j)\big)\, \Big)
\\
&\,=\, 3(1-w^2)\,(a_1 + a_2)
\\
&
\\
&\quad\,+\, \frac  12 \Big[ F''( 1  - s_2 a_2  ) -F''( w- s_3(a_1+a_2 ) ) \Big]\,(a_1^2 + a_2^2)
\, + \, O\left (  \ve^{2+\tau} e^{-\sigma|t|} \, \right)
\\
&
\\
&\, =\, 3( 1- w^2 )\, (a_1 + a_2)
 \,+\,  O\left (  \ve^{2+\tau} e^{-\sigma|t|} \, \right) \, .
\end{align*}
Hence, recalling relations \equ{fj}, \equ{rho}, the definitions of $a_1$, $a_2$ and the asymptotic expansions \equ{r5}, \equ{r6} for $j=\ell-1$ and $j=\ell+1$,
we find

$$
F(u_0) - \sum_{j=1}^m F(w_j(t - f_j))\  =
$$

\be  6(-1)^{\ell-1}(1- w^2(t))\ve^2 \rho_\ve \,\left  [ \, e^{- \sqrt{2}(h_\ell  -h_{\ell -1})} e^{ \sqrt{2} t}\, -
 e^{ -\sqrt{2}(h_{\ell +1}- h_\ell ) } e^{ -\sqrt{2} t}\, \right ] \, +\, \theta_\ell  .
\label{ee} \ee

\noindent
where $ \theta_\ell  =  O\left (  \ve^{2+\tau} e^{-\sigma|t|} \, \right) $.

\medskip
Substituting \equ{ee} in expression \equ{expr1} we then find
\begin{align}
&(-1)^{\ell-1}S(u_0) (y,t+ f_\ell)\nonumber
\\
&\,=\,
6(1- w^2(t))\ve^2 \rho_\ve \,\left[ \, e^{- \sqrt{2}(h_\ell  -h_{\ell -1})} e^{ \sqrt{2} t}\,-\,
e^{ -\sqrt{2}(h_{\ell +1}- h_\ell ) } e^{ -\sqrt{2} t}\, \right ] \   \nonumber
\\
&\quad\,-\, \ve^2\, \Bigl(\Delta_{\KK} h_\ell  +  (t + f_\ell  )K\,\Bigr) w'(t)
+\, \ve^2|\bigtriangledown_{\KK }h_\ell  |^2 w''(t)\nonumber
\\
&\quad\,+\, (-1)^{\ell-1}\Theta_\ell (\ve y, t)\,.
\label{expr2-1}
\end{align}
Here we have denoted

\begin{align}
\Theta_\ell (\ve y, t)\,=\ &  \theta_\ell (\ve y, t)-\, \ve^2\sum_{j\ne \ell} \,\Bigl(\Delta_{\KK} h_j +(t+ f_\ell )K\,\Bigr) w_j'(t + f_\ell  -f_j)
\nonumber
\\
&\,+\,\ve^2\sum_{j\ne \ell}|\bigtriangledown_{\KK }h_j |^2 w_{j}''(z- f_j)\nonumber
\\
&\,+\,\ve^3 (-1)^{\ell-1}z\, a_{ik}^1\,\pp_ih_j \,\pp_k h_j \, w_j''(t+ f_\ell -f_j)\nonumber
\\
&\,+\,\ve^3(t+ f_\ell )\, \big[\, a^1_{ik}  \,\pp_{ik} h_j   \,+ \,  b_i^1\,\pp_i h_j \,\big]\, w_j'(t + f_\ell  -f_j)
\nonumber
\\
&\,+\,\ve^3 \big[\, z^3\,a^3_{NN}\,  +   z^2\, a_{iN}^2 \, \pp_{i} h_j\big] \,  w_j''(t + f_\ell  -f_j),
\label{expr3}
\end{align}

\noindent
where the coefficients are understood to be evaluated at $\ve y$ or $\big(\ve y,\ve (t+ f_\ell (\ve y))\big)$.

\medskip
While this expression has been obtained assuming $2\le \ell\le m-1$, we see that it also holds for $\ell= m$, $\ell=1$. The cases $\ell=1$ and $\ell=m$ are dealt similarly. The only difference is that the term  $\left  [ \, e^{ -\sqrt{2}(h_\ell  -h_{\ell -1})} e^{ \sqrt{2} t}\, -
 e^{ -\sqrt{2}(h_{\ell +1} -h_{\ell }) } e^{ -\sqrt{2} t}\, \right ]$ gets respectively replaced by
 \be
   -
 e^{ -\sqrt{2}h_{1} } e^{ -\sqrt{2} t}\, ,\, \ell=1
\quad \hbox{ and }\quad \, e^{ -\sqrt{2}h_{m}} e^{ \sqrt{2} t}\, ,\, \ell =m .
\label{esp}\ee

\medskip
\subsection{Size of the error}
Examining expression \equ{expr3} we see that the error in the considered region is made up by terms that can be bounded by a power
of $\ve$ times a factor with  exponential decay in $t$. We introduce the following norm
for a function $g(y,t)$  defined on $\KK_\ve\times\R$.
Let $\sigma>0$, $1<p\le +\infty$. We set

\be
\|g\|_{p,\sigma} = \sup_{(y,t) \in\KK_\ve\times\R}  e^{\sigma|t|} \|g\|_{L^p\big( B( (y,t), 1)\big)}.
\label{norm}
\ee

We want to consider the error associated to points in the set $A_\ell $ as a function defined in
 the entire space $\KK_\ve\times\R$. To do so, we consider a smooth cut-off function
 $\zeta (s) $ with $\zeta(s) =1$ for $s<1$ and   $\zeta(s) =0$ for $s>2$ and define
$$
\zeta_\ve (t) = \zeta \big( |t| - \frac {\rho_\ve} 2 - 2M \big).
$$
We extend the error as follows.
Let us set
\begin{align}
S_\ell (u_0) \,:=\ &
6\big(1- w^2(t)\big)\ve^2 \rho_\ve \,\left  [ \, e^{ -\sqrt{2}(h_\ell  -h_{\ell -1})} e^{ \sqrt{2} t}\, -
e^{ -\sqrt{2}(h_{\ell +1} -h_{\ell }) } e^{ -\sqrt{2} t}\, \right ]\, \zeta_\ve(t)\nonumber
\\
&\,-\, \ve^2\, \Bigl(\Delta_{\KK} h_\ell  +  (t + f_\ell  )K\,\Bigr) w'(t)
\,+\, \ve^2|\bigtriangledown_{\KK_\ve }h_\ell  |^2 w''(t) \nonumber
\\
&\,+\,(-1)^{\ell-1}\,\zeta_\ve(t) \,\Theta_\ell (\ve y, t)\,,
\label{Sn}
\end{align}
where the cut-off expressions are understood to be zero outside the support of $\zeta_\ve$. We see that
$$
(-1)^{\ell-1}S(u_0) (\ve y,t+ f_\ell)\, =   S_\ell (u_0) (y,t)\foral (y,t)\in A_\ell .
$$

The following lemma on the accuracy of the error is readily checked.
\begin{lemma}\label{error1}
For a given $0<\sigma< \sqrt{2}$ and any $p>1$ we  have the  estimates
\begin{equation}
{\|}S_\ell (u_0 ){\|}_{p,\sigma}\leq C\, \ve^{2-\tau} ,
\quad
{\|} \zeta_\ve \,\Theta_\ell  {\|}_{p,\sigma}\leq C\, \ve^{3-\tau},
\label{lerr2}
\end{equation}
where $\tau$ is any number with $\tau > \frac 1{2\sqrt{2}} \sigma $ and $\tau > \frac {N-1}p$.
\end{lemma}
\proof
The proof amounts to a straightforward verification of the bound term by term.
Let us consider for instance
$$ E_1 =  6\big(1- w^2(t)\big)\ve^2 \rho_\ve \,\left  [ \, e^{ -\sqrt{2}(h_{\ell } -h_{\ell -1})} e^{ \sqrt{2} t}\, -
 e^{ -\sqrt{2}(h_{\ell +1} -h_{\ell })} e^{ -\sqrt{2} t}\, \right ]\, \zeta_\ve(t) .$$
Then for $|t|\le \frac{\rho_\ve}2$ we get
\begin{align*}
|E_1|
\,\le\,
  C \ve^2|\log\ve|
\,\le \,
C e^{- \sqrt{2}{\rho_\ve}}
\, \le \,
e^{- \sigma \frac {\rho_\ve} 2}\, e^{ -(\sqrt{2} - \frac 12 \sigma)\rho_\ve}
\, \le\,
 e^{- \sigma |t| }\, \ve^{2-\tau},
\end{align*}
where  $\tau >  \frac 1{2\sqrt{2}} \sigma $. This implies $\|E_1\|_{p,\sigma} \le C\ve^{2-\tau}$ for any $1<p\le +\infty$. Now,
let us consider the term
$$
E_2(y,t)=
\ve^2\, \Delta_{\KK} h_\ell (\ve y)\, w'(t).
$$
Then for any $\sigma \le \sqrt{2}$ we have
\begin{align*}
e^{\sigma |t|} \|E_2\|_{L^p( B((t,y),1)}  &\,\le\,  C \ve^2 \|\Delta_{\KK} h_\ell (\ve \ttt y)\|_{L^p( B(y,1)}
\\
&\,\le\,C \ve^{ 2-\frac {N-1} p} \|\Delta_{\KK} h_\ell \|_{L^p( \KK )}
\, \le \, C \ve^{ 2-\frac {N-1} p}.
\end{align*}
The rest of the terms are dealt similarly, being in fact roughly at least $\ve$ times  smaller than those above. \qed

\medskip
Very important for subsequent developments is the Lipschitz character of the error in the parameter function
$\hh = (h_1,\ldots, h_N)$.
Let us write $S_j(\hh)$ to emphasize the dependence on this function.
We have

\begin{lemma}\label{lispchitzerror}
Let us assume that the vector-valued functions $\hh_1$, $\hh_2$ satisfy the constraints \equ{assh}. We have the validity of the following Lipschitz conditions.
$$
\|  S_j(\hh^1) - S_j(\hh^2)\|_{p,\sigma} \le   C\ve^{2-\tau} \|\hh^1-\hh^2\|_{W^{2,p}(\KK)},
$$
$$
\|\zeta_\ve \,\Theta_\ell (\hh^1)  - \zeta_\ve \,\Theta_\ell (\hh^2) \|_{p,\sigma} \le  C\ve^{3-\tau} \|\hh^1-\hh^2\|_{W^{2,p}(\KK)},
$$
for $\tau > \frac 1{2\sqrt{2}} \sigma $ and $\tau > \frac {N-1}p$.

\end{lemma}

\proof
Again the proof consists in establishing the bound for each of its individual terms, more precisely, we need to bound now
for instance
$ \pp_{\pp_i h_j}S_j(\hh)$. Since the dependence on this object, and as well on second derivatives comes in linear or quadratic way,
always multiplied by exponentially decaying factors and small powers of $\ve$, the desired result directly follows. The dependence on
the values of the functions $h_j$ appears in a more nonlinear fashion, however smooth and exponentially decaying.
We omit the details. The complete arguments are rather similar to those in the proof of Corollary 5.1 of \cite{dkpw2}. \qed

\medskip
\subsection{The global approximation}

The approximation $u_0$  is so far defined only in a neighborhood of $\KK_\ve$ in $\MM_\ve$, where  the local Fermi coordinates make sense.
Let us assume that $m$ is an odd number. In this case we require that  $\KK_\ve$ separates $\MM_\ve$ into two components that we denote
$\MM_\ve^-$ and $\MM_\ve^+$.

 Let us use the convention that the normal to $\KK_\ve$  points in the direction of $\MM_\ve^+$. Let us consider the function
 $\HH$ defined in $ \MM_\ve \setminus  \KK_\ve$ as
\be
\HH(x)\, := \,
\left\{
\begin{matrix}
\ 1 &\quad\hbox{if } x\in \MM_\ve^+,
\\
-1 &\quad\hbox{if } x\in  \MM_\ve^-.
\end{matrix}
\right.
\label{HH}\ee
Then our approximation $u_0(x)$ approaches $\HH(x)$ at an exponential rate $O( e^{-\sqrt{2} |t|})$ as $|t|$ increases.
 The global approximation we will use  consists simply of interpolating $u_0$ with $\HH$ sufficiently well-inside $ \MM_\ve \setminus  \KK_\ve$  through a cut-off in $|z|$.  Let
 $\NN_\delta$ be the set of points in $\MM_\ve$ that have Fermi coordinates $(y,z)$ well-defined and $|z| < \frac \delta\ve$.
with some positive constant $\delta<\delta_0/10$.

\medskip
Let $\eta(s)$ be a smooth cut-off function with $\eta(s) =1$ for $s<1$ and $=0$ for $s>2$
and define
\be
\eta_\delta (x) \, := \,
\left\{
\begin{matrix}
\ \eta( \, |z| - \frac \delta\ve )  &  \quad\hbox{if } x\in \NN_\delta \, , \\
 0 &\quad\hbox{if } x \not \in  \NN_\delta.
\end{matrix}
\right.
\label{etadelta}\ee
 Then we let our global approximation $\ww(x)$ be simply defined as
\be
\ww  \, := \, \eta_\delta u_0  \,+\, (1-\eta_\delta )\HH,
\label{global}\ee
where  $\HH$ is given by \equ{HH} and $\ww$ is just understood to be $\HH(x)$  outside $\NN_\delta$.

\medskip
Since $\HH$ is an exact solution in $\R^N\setminus \MM_\ve$, the global error of approximation is simply computed as

\be
S(\ww) \,=\, \Delta \ww \,+\, F(\ww)
\,=\, \eta_\delta  S(u_0)\,+\, E,
\label{error}\ee
where
$$
E\,=\, 2\nabla \eta_\delta\nabla u_0
\,+\, \Delta \eta_\delta (u_0 -\HH)
\,+\,F\big(\,\eta_\delta u_0  + (1-\eta_\delta )\HH\,\big)
\,-\, \eta_\delta F(u_0) \, .
$$
Observe that $E$ has exponential size  $O(e^{-\frac c \ve})$ inside its support, and hence the contribution of this error to the entire error
is essentially negligible.

\medskip
If $m$ is even, we simply define
\be
\ww  \, := \, \eta_\delta u_0  + (1-\eta_\delta )(-1).
\label{global2}\ee
In this case there is no need that $\KK$ separates $\MM$ into two components.

\setcounter{equation}{0}
\section{The gluing procedure}\label{section3}
\setcounter{equation}{0}

Once the global approximation $\ww(x)$ in \equ{global} or \equ{global2} has been built, we then want to find
a solution to the full problem of the form
$$
u(x) = \ww(x) +\vp(x),\quad x\in \MM_\ve
$$
where $\vp(x)$ is a small function.
Thus $\vp$ must satisfy
\be
\Delta_{\MM_\ve} \vp + F'(\ww)\vp = -S(\ww ) - N(\vp) \inn \MM_\ve
\label{p111}\ee
where
$$
N(\vp) = F(\ww +\vp) -F(\ww) - F'(\ww)\vp.
$$
We shall look for a solution of the form
$$
u(x) \,=\, \sum_{j=1}^m \zeta_{j2}(x)  \ttt\phi_j(y,z ) \,+\, \psi(x),
$$
where the functions $\ttt \phi_j$ are defined in the entire space $\KK_\ve\times\R$.
Then the equation is equivalent to
\begin{align*}
&\sum_{j=1}^m \zeta_{j2}\Big[\Delta_{\MM_\ve} \ttt\phi_j\,+\, F'(\ww)\ttt\phi_j \,+\, \zeta_{j1} (F'(\ww) +2) \psi\,+\, \zeta_{j1}N(\psi + \ttt\phi_j)\, +\, S(\ww ) \,\Big]
\\
&\,+\,\sum_{j=1}^m\Big[ 2\big<\nn_{\MM_\ve}  \zeta_{j2},\nn_{\MM_\ve} \ttt\phi_j\big>  +  \ttt\phi_j \Delta_{\MM_\ve} \zeta_{j2}\Big]
\,+\, \Delta \psi \,-\, \Big(2-F'(\ww)\big(1- \sum_{j=1}^m\zeta_{j1}\big)\Big) \psi
\\
&\,+\,\big(1- \sum_{j=1}^m\zeta_{j1}\big) \, N\Big( \psi + \sum_{i=1}^m\zeta_{i2} \ttt\phi_i\Big)
\,+\,\Big(1- \sum_{j=1}^m\zeta_{j2}\Big)\, S(\ww) \,=\, 0 \inn \MM_\ve.
\end{align*}
 This system will be satisfied if the $(m+1)$-tuple $(\ttt\phi_1,\ldots, \ttt\phi_m,\psi)$ solves the system
\be
  \Delta_{\MM_\ve} \ttt\phi_j\,+\, F'(\ww)\ttt\phi_j
  \,+\,  \zeta_{j1} \big(F'(\ww) +2\big) \psi
  \,+\, \zeta_{j1}N\big(\psi + \ttt\phi_j\big)
  \,+\,  S(\ww ) \,=\, 0,
  \label{glu1}
\ee
for $ |z|< C|\log\ve|$, $j=1,\ldots, m$, and
\be
  -\Delta_{\MM_\ve} \psi \, +\, 2\psi \, = \, \QQ ( \psi, x) \inn \MM_\ve,
  \label{glu2}
\ee
where we have denoted
\begin{align}
\QQ(\psi,x) \, :=\,&
\Big(1- \sum_{j=1}^m \zeta_{j1}\Big) \left \{ \,\big[ 2 +F'(\ww)\big]\, \psi
\,+\,N\big( \psi + \sum_{i=1}^m \zeta_{i2} \ttt\phi_i\big)  \, \right \}
\label{glu3}
\\
&
\,+\,\big(1- \sum_{j=1}^m \zeta_{j2}\big) S(\ww)
\,+\,\sum_{j=1}^m\Big[\, 2\big<\nn_{\MM_\ve}  \zeta_{j2},\nn_{\MM_\ve} \ttt\phi_j\big>  +  \ttt\phi_j \Delta_{\MM_\ve} \zeta_{j2} \Big] \, .
\nonumber
\end{align}

The gluing procedure consists in solving equation \equ{glu2} for $\psi$
in terms of a given ${\ttt \phi}=(\ttt\phi_1,\ldots, \ttt\phi_m) $ chosen arbitrary but sufficiently small,
and then substituting the result in equation \equ{glu1}.
Let us assume the following constraints on the $\ttt\phi_j$'s:
\be
\ttt \phi_j(y,z)\ =:\ \phi_j\big(y, z- f_j(\ve y)\big) ,
\quad
\|\phi_j\|_{2,p,\sigma} \le 1\foral j=1,\ldots m.
\label{41}
\ee

\begin{lemma}\label{lemapsi}
Given functions $\phi_j$ and $\hh$ satisfying respectively constraints \equ{41} and \equ{assh}, there exists a unique solution $\psi = \Psi(\phi,\hh)$ to equation
\equ{glu2} with
$$
\|\psi\|_\infty \le C(\ve^{4-\tau}+  \ve^{2-\tau}\|\phi\|_{2,p,\sigma}),
$$
for a small $\tau>0$. In addition the operator $\Psi$ satisfies the Lipschitz condition
\be
\|\Psi(\phi^1,\hh^1)-\Psi(\phi^2,\hh^2)\|_\infty
\,\le\,
C \ve^{2-\tau}\big[\, \|\phi^1 - \phi^2\|_{2,p,\sigma}  +  \|\hh_1 -\hh_2\|_{2,p}\, \big].
\label{lippsi}
\ee

\end{lemma}

\proof
Let us consider first the linear equation
\be
 -\Delta_{\MM_\ve}  \psi   + 2\psi  = E(x) \inn \MM_\ve.
\label{eqpsi}\ee
We claim that if we set
$$
\|E\|_{p,0} = \sup_{x\in M_\ve } \| E\|_{L^p( B(x,1))},
$$
then problem \equ{eqpsi} has, for all small $\ve>0$, a unique bounded solution $\psi = \mathfrak{A}(E)$, which in addition satisfies
$$
\|D_{\MM_\ve}\psi\| + \|\psi\|_\infty
\,\le\,
C \|E\|_{p,0},
$$
provided that $p>m$. To prove this claim, it suffices to establish the a priori estimate in $L^\infty$-norm. If that was not true, there would be sequences $\ve=\ve_n$,
$\psi_n$, $E_n$, with $\|E_n\|_{p,0}\to 0$, $\|\psi_n\|_\infty =1$ such that
$$
 -\Delta_{\MM_\ve}  \psi_n   + 2\psi_n  = E_n\inn \MM_\ve.
$$
Using local normal coordinates around a point $p_n\in \MM_\ve$ where $|\psi_n(p_n)|=1$, the same procedure as in the proof of
the a priori estimate in Proposition \ref{prop1}  leads us to local convergence of $\psi_n$ to a nontrivial bounded solution of
$$
-\Delta_{\R^N} \psi +2 \psi = 0 \inn \R^N,
$$
and a contradiction is reached.

\medskip
To solve equation \equ{glu2} we write it in fixed point form as
\be
\psi = \mathfrak{A} ( \QQ (\psi, \cdot)) .
\label{fp1}\ee
In the region where the functions $(1-\sum_i \zeta_{1i})\zeta_{2j}$,   $\nn_{\MM_\ve}\zeta_{2j}$, $\Delta_{\MM_\ve}\zeta_{2j}$ are supported
we have, thanks to \equ{41},
$$
|\ttt \phi (x)| + |\nn_{\MM_\ve} \ttt \phi(x)|
\,\le\,
C\,e^{-\sigma{ |\xi_j -\xi_{j-1}|} }\,\|\phi_j\|_{2,p,\sigma}
\, \le\,
\ve^{2-\tau}\|\phi_j\|_{2,p,\sigma},
$$
for a  small $\tau >0$.
We also notice that
$$
|\QQ(0,x)| \le  C\ve^{4-\tau}\Big(|D^2_\KK \hh(\ve y)| + |D_\KK \hh(\ve y)| + |\hh(\ve y)|+ 1\Big)
\,+\,
\ve^{2-\tau}\sum_{j=1}^m\|\phi_j\|_{2,p,\sigma}.
$$
We observe then that
$$
\|\QQ(0,\cdot)\|_{0,p} \le  C\ve^{4-\tau - \frac Np} \|h\|_{W^{2,p}(\KK)}   +  \ve^{2-\tau}\|\phi\|_{2,p,\sigma}\le C(\ve^{4-\tau}+  \ve^{2-\tau}\|\phi\|_{2,p,\sigma})\ .
$$
We check next the Lipschitz character of this operator, not just in $\psi$ , but also in the rest of its arguments.
Let us write $\QQ = \QQ(\psi, \hh, \phi)$ and assume
\be
\|\phi\|_{2,p,\sigma} \le 1, \quad \|\psi\|_{2,p,\sigma} \le \beta \ve^{2-\tau},\quad \|\hh\|_{W^{2,p}(\KK)} \le M.
\label{ass}\ee
We consider  $(\psi^l, \phi^l,\hh^l )$, $l=1,2$, satisfying \equ{ass},  and denote $\QQ^l = \QQ(\psi^l,\phi^l, \hh^l)$.
We will show that
\begin{align}
&\|\QQ(\psi^1,\phi^1, \hh^1)- \QQ(\psi^2,\phi^2, \hh^2)\|_{0,p}\nonumber
\\
&\quad\,\le\, C \ve^{2-\tau}\big[\,\|\psi_1 - \psi_2\|_{\infty}  + \|\phi^1 - \phi^2\|_{2,p,\sigma} + \|\hh_1-\hh_1\|_{2,p}\,\big].
\label{lipQQ}
\end{align}
Let us observe that  for $(\psi, \phi,\hh)$ satisfying \equ{ass},
$$\QQ(x) = \QQ (\psi, \phi,\hh)(x)\,=\,
  \QQ\Big(\psi(x) ,\hh(x), D\hh(x) , D^2\hh(x), \phi(x), D\phi(x) , x\Big).
$$
We decompose
$$\QQ(x) \ =\   \underbrace { \Big(1- \sum_{j=1}^m \zeta_{j2}\Big) S(\ww) }_{Q_1}
\,+\,
\underbrace{\sum_{j=1}^m\Big[\, 2\big<\nn_{\MM_\ve}  \zeta_{j2},\nn_{\MM_\ve} \ttt\phi_j\big>
+  \ttt\phi_j \Delta_{\MM_\ve} \zeta_{j2} \Big]}_{Q_2}
$$

\be \underbrace{
  \,+\,\Big(1- \sum_{j=1}^m \zeta_{j1}\Big)\,\left( \,\big[ 2 +F'(\ww)\big]\, \psi\,+\,
     N\Big( \psi + \sum_{i=1}^m \zeta_{i2} \ttt\phi_i\Big) \right) }_{Q_3} \ .
 \label{glu8}\ee
Then we find
$$
\pp_{\psi} \QQ  =  \pp_{\psi}\QQ_3 = \Big(1- \sum_{j=1}^m \zeta_{j1}\Big)
\left \{ \,\big[ 2 +F'(\ww)\big]\,+\,N'\Big( \psi + \sum_{i=1}^m \zeta_{i2} \ttt\phi_i\Big) \, \right \}\
$$
where $N'(s) = F'(\ww + s ) - F'(\ww) = O(|s|)$.
It follows that $\pp_{\psi} \QQ  = O(\ve^{2-\tau})$ in the considered range for the parameters.
Now,
$$
\pp_{h_k} \QQ_3  =   \Big(1- \sum_{j=1}^m \zeta_{j1}\Big) \left \{ \,F''(\ww) \, \psi
\,+\,N'\Big( \psi + \sum_{i=1}^m \zeta_{i2} \ttt\phi_i\Big)  \, \right \}\, w'(z- f_j(\ve y)),
$$
since  $\pp_{h_k} \ww = w'\big(z- f_j(\ve y)\big) $.
Thus  $\pp_{h_k} \QQ_3=  O(\ve^{4-\tau})$.
We also have  $\pp_{h_k}\QQ_1 = O(\ve^{4-\tau})$ and
$$
\pp_{ D^l h_k} \QQ  =  \pp_{ D^l h_k} \QQ_1 =   \Big(1- \sum_{j=1}^m \zeta_{j2}\Big)  \pp_{ D^lh_k}S(\ww) = O(\ve^{4-\tau}),  \quad l=1,2.
$$
Finally,
$$
\pp_{ \phi_j } \QQ   =  \pp_{ \phi_j } \QQ_3 =   \Delta_{\MM_\ve} \zeta_{j2} =O(1) ,
\quad
\pp_{ D \phi_j } \QQ  =  2\nn \zeta_{j2} = O(1).
$$
As a conclusion, using  the mean value formula and the facts
\begin{align*}
|\phi^1 -\phi^2| +  |D\phi^1 -D\phi^2| \le e^{-\sigma|t|} \|\phi^1 -\phi^2\|_{2,p,\sigma},
\\
\\
\| D^2\hh^1(\ve y)  -D^2\hh^2 (\ve y)\|_{0,p} \le  C\ve^{-\frac mp} \| D^2\hh^1 -D^2\hh^2\|_{L^p(\KK)},
\end{align*}

\noindent we readily find the validity of \equ{lipQQ}.
In particular, we obtain
that for $\|\psi_l\|_\infty \le  \beta \ve^{2-\tau}$, $l=1,2$, and
\be
\| \QQ (\psi_1, \phi,\hh ) -\QQ(\psi_2,\phi,\hh)\|_{p,0} \le C\ve^{2-\tau}   \|\psi_1-\psi_2\|_\infty.
\label{lipQ}\ee
Thus, from the contraction mapping principle, we find that for certain $\beta>0$ large and fixed, problem \equ{fp1} has a unique solution
$\psi = \Psi(\phi,\hh)$ such that
\be
 \|\Psi(\phi,\hh)\|_\infty \le C(\ve^3+  \ve^{2-\tau}\|\phi\|_{2,p,\sigma}).
\label{psi1}\ee
The Lipschitz dependence of $\Psi$ \equ{lippsi} in its arguments  follows  immediately from \equ{lipQ} and the fixed point characterization \equ{fp1}. \qed

\bigskip
Now, assuming that $\|\phi\|_{2,p,\sigma}$ is in the considered range, we substitute $\psi = \Psi(\phi,\hh)$ in \equ{glu1} and then obtain that
$$
\vp = \Psi(\phi,\hh) + \sum_{j=1}^m \zeta_{j2}\ttt \phi_j ,
\qquad
\ttt\phi_j(y,z) = \phi_j(y, z-f_j(\ve y)),
$$
solves problem \equ{p111} if  the vector $\phi = (\phi_1,\ldots,\phi_j)$ satisfies the system of equations
\be
\Delta_{\MM_\ve} \ttt\phi_j + F'(\ww)\ttt\phi_j +  \zeta_{j1}\big(F'(\ww) +2\big) \psi
+ \zeta_{j1}N\big(\Psi(\phi,\hh) + \ttt\phi_j\big) + \zeta_{j1} S(\ww) = 0,
\label{glu11}
\ee
in the support of $\zeta_{j2}$.
We want to extend these equations to the entire $\KK_\ve\times\R$.
We recall that in $(y,z)$ coordinates we can write
$$
\Delta_{\MM_\ve} = \partial_{zz}^2 + \Delta_{\KK_\ve}  - \ve^2zK(\ve y)\pp_z +  B,
$$
where $B$ is a small operator given by \equ{defB}.
It is convenient to rewrite  equations \equ{glu11} in terms of the functions $\phi_j$ defined as
$$
\phi_j(y,t) = \ttt \phi_j(y,t+ f_j(\ve y)).
$$
We find in coordinates $(y,t)$,

\be
\Delta_{\MM_\ve} \ttt\phi_j  \,=\, \partial_{tt}^2\phi_j
\,+\, \Delta_{\KK_\ve}\phi_j
\,+\, B_j^1\phi_j
\,+\, B_j^2\phi_j,
\label{lap1}
\ee
where
\begin{align*}
B_j^1\phi  \ =\ &\ve^2\big|\nn_\KK h_j(\ve y)\big|^2\, \partial_{tt}^2\phi
\,-\,\ve^2\Delta_\KK h_j(\ve y)\, \pp_t\phi
\\
&\,-\,\ve^2 K(\ve y)\, \big(t+ f_j(\ve y)\big) \pp_t\phi
\,-\, 2\ve \big< \nn_{\KK} h_j(\ve y)\,,\,  \nn_{\KK_\ve} \pp_t\phi  \big>,
\end{align*}
and, expressed in local coordinates $(\py,t)$,  $y= Y_p(\py)$, the operator $B_j^2\phi$ becomes

\begin{align}
B_j^2\phi \,= \ &\ve^3 \big(t+ f_j(\ve \py) \big)^3\,a^3_{NN}\,\pp_{tt}
\,+\, \ve^2\, \big(t+ f_j(\ve \py)\, \big) \, b_i^1 \, \big( \pp_i \phi -  \ve \pp_i h_j(\ve \py)\, \pp_t \phi\,\big)
\nonumber
\\
\nonumber
\\
&\,+\, \ve \,\big(t+ f_j(\ve \py)\big)\, a_{il}^1\,
\Big[ \, \pp_{il} \phi  -2\ve \pp_l h_j(\ve \py)\, \pp_{it} \phi
\label{Bj}
\\
\nonumber
\\
&\qquad\qquad\qquad\qquad\qquad- \ve^2 \pp_{il} h_j(\ve \py) \, \pp_t \phi\, +\, \ve^2\pp_ih_j\pp_l h_j (\ve\py)\, \pp_{tt}\phi \,\Big]
\nonumber
\\
\nonumber
\\
 &\, +\,\ve^2 \big(t+ f_j(\ve \py) \big)^2\,a^2_{iN}\,\big( \pp_{it} \phi -  \ve \pp_i h_j(\ve \py)\, \pp_{tt} \phi\,\big)
 \,+\, \ve^3 \,\big(t+ f_j(\ve \py) \big)^2\, b_N^2\,\pp_t\,  ,
\nonumber
\end{align}

\noindent
with coefficients evaluated at $(\ve \py, \ve t+ \ve f_j(\ve\py))$.
The difference between the operators $B_j^1$ and $B_j^2$ is that the expression for $B_j^1$ actually makes sense for all $(y,t)$, while $B_j^2$ does only up to
$|t| <\delta/\ve$. We set $\chi_0(t) = \zeta(  |t| -  10 \log\ve)$,
where, we recall, $\zeta(\tau ) = 1$ for $\tau <1$ and $=0$ for $\tau >2$. Then we extend the operator $B_j^1 + B_j^2$ to entire space $(y,z)$
setting  $$B_j := B_j^1 + \chi_0 B_j^2. $$

Let us relabel

$$
\chi_{js} (y,t)
\,:=\, \zeta_{js} (t+ f_j(\ve y))
\,=\, \zeta ( | t+ h_j(\ve y)| -  d_\ve  -  s )  ,
$$
and denote

$$
\Psi_j(\phi,\hh) (y,t)\,:=\, \Psi(\phi,\hh)( y, t+f_j(\ve y)) ,  \quad
$$

\be
S_j(\hh )(y,t)  \,:=\,   \chi_{j3}S(\ww)( y, t+ f_j) , \
\label{Sjj}
\ee
(observe that this is the same $S_j$ introduced in \equ{Sn})

$$
\ww_j :=  \ww(y, t+f_j) = w(t) + \theta_j,
$$

\noindent
where $\theta_j(y,t) =  O(\ve^{-2+\tau})$.
We have

\begin{align}
&\partial_{zz}^2\ttt\phi_j + \Delta_{\KK_\ve}\ttt\phi_j +  F'(w_j)\ttt\phi_j
\,+\,\zeta_{j3}\, B\ttt\phi_j  +  \zeta_{j3}\,\big(F'(\ww) - F'(w_j)\big)\ttt\phi_j  \nonumber
\\
&+\zeta_{j1}\,\left[\,\big(F'(\ww) +2\big) \Psi(\phi,\hh)  +  N\big(\Psi(\phi,\hh) + \ttt\phi_j\big) +  S(\ww )\right ] = 0 \inn \KK_\ve\times\R ,
  \label{glu111}
\end{align}

\noindent
where $w_j(y,z) = w(z-f_j(\ve y))$.
Finally, we recast equations \equ{glu111} as
\be
  \partial_{tt}^2\phi_j + \Delta_{\KK_\ve}\phi_j +  F'(w(t))\ttt\phi_j  +  S_j(\hh ) + \, \NNN_j(\phi,\hh) \, = 0 \inn \KK_\ve\times\R,
  \label{glufin}
\ee
for all $j=1,\ldots, m,$ where
\be
  \NNN_j(\phi,\hh)\: =\  \BB_j(\phi_j)
  \,+ \, \chi_{j1}\, \Big[  \,\big(F'(\ww_j) +2\big) \Psi_j(\phi,\hh)+  N\big(\Psi_j(\phi,\hh) + \phi_j\big) \,\Big],
  \label{Nj}
\ee
  with
\be
\BB_j(\phi_j) = \chi_{j3}\,\Big[  B_j \phi_j  + \big(F'(\ww_j) - F'(w) \big)\phi_j \Big],
\quad
B_j\hbox{ given by \equ{Bj}}.
\label{BBj}
\ee
We will concentrate in what follows in solving system \equ{glufin}.
We will do this in two steps: 1. solving a projected version of the problem, carrying $\hh$ as a parameter,
and 2. finding $\hh$ such that the solution of this projected problem
is an actual solution of \equ{glufin}.
We consider then the system, for all $ j=1,\ldots, m$

\begin{align}
\begin{aligned}
\partial_{tt}^2\phi_j\, + \, \Delta_{\KK_\ve }\phi_j\,& = \, -S_j (\hh) - \NNN_j(\phi,\hh) \,  +\, c_j (y) w'(t)
\quad \hbox{ in } \KK_\ve \times \R
,\\ \\
\int_\R \phi_j(y,t)\,w'(t)\,\mathrm{d}t&=0\quad\mbox{on } \KK_\ve \, ,
\quad
c_j(y) =  \frac {\int_\R \big( S_j (\hh ) + \NNN_j(\phi,\hh) \big)\, w'\,\mathrm{d}t }{\int_\R {w'}^2\mathrm{d}t } \, .
\end{aligned}
\label{glufinp}
\end{align}

To solve it we need a suitable invertibility theory for the linear operator involved in the above equation. We do this
next.

\setcounter{equation}{0}
\section{The auxiliary linear projected problem}

Crucial for later purposes is a solvability theory for the following linear problem:

\begin{align}
\begin{aligned}
\partial_{tt}^2\phi\, + \, \Delta_{\KK_\ve }\phi\,& = \, g(y,t)\,  +\, c (y) w'(t)
\quad \hbox{ in } \KK_\ve \times \R
,\\ \\
\int_\R \phi(y,t)\,w'(t)\,\mathrm{d}t&=0  \foral  y\in  \KK_\ve \, ,
\quad
c(y) =  - \frac {\int_\R g(y,t ) w'\mathrm{d}t }{\int_\R {w'}^2\mathrm{d}t } \, .
\end{aligned}
\label{p1}
\end{align}

We have the following result.

\begin{prop}\label{prop1}
Given $p>m$  and $0< \sigma < \sqrt{2}$, there exists a constant
$C>0$ such that for all sufficiently small $\ve  >0$ the following
holds. Given $g$ with $\|g\|_{p ,\sigma}< +\infty$, then Problem
$\equ{p1}$ has a unique  solution $\phi$ with $\|\phi
\|_{\infty,\sigma} < +\infty$, which  in addition satisfies \be
\| D^2 \phi\|_{p,\sigma } + \|D\phi\|_{\infty , \sigma} + \|\phi\|_{\infty , \sigma }
\, \le \, C \|g\|_{p ,\sigma}\, .
\label{cota} \ee
\end{prop}

The main fact needed is that the one-variable solution
$w$ of \equ{definitionofH} is {\em nondegenerate} in  $L^\infty(\R^{m})$ in the sense that
 the linearized operator
$$
L(\phi) = \Delta_y\phi + \partial_{tt}^2\phi + F'(w(t))\phi ,\quad (y,t)\in\R^{N-1}\times\R,
$$
is such that the following property holds.
\begin{lemma}\label{lemma l1}
Let $\phi$ be a bounded, smooth solution of the problem
\begin{equation}
L(\phi) = 0 \quad\hbox{in } \R^{N-1}\times \R.
\label{l3}\end{equation}
Then $\phi(y,t) = Cw'(t) $ for some $C\in \R$.
\end{lemma}

\proof
This fact is by now standard, so we only sketch the proof.
The  one-dimensional operator
 $
L_0(\psi)= \psi'' + F'(w)\psi
$
is such that $L_0(w') =0$ and $w'>0$, hence $0$ is its least eigenvalue. Using this, it is easy to show that
there is a constant $\gamma >0$ such that
whenever $\int_\R \psi w' = 0$ with $\psi\in H^1(\R)$ we have that
\be
\int_\R \big(\, |\psi'|^2 - F'(w)\psi^2\,\big)\, \mathrm{d}t\, \ge \, \gamma
\int_\R (\, |\psi'|^2 + |\psi|^2\,)\,\mathrm{d}t . \label{q}
\ee
Let
$\phi$ be a bounded solution of equation \equ{l3}.
Since $F'(w(t))\sim -2$ for all large $|t|$ an application of the maximum principle shows that if $0<\sigma< \sqrt{2}$ and $t_0>0$ is large then

$$
|\phi(y,t)| \ \le\ C\|\phi\|_\infty   e^{-\sigma|t|}  \quad \hbox{if } |t| >t_0 \, .
$$
On the other hand, the function
$$
{\bar\phi}(y,t) \,=\,
\phi(y,t)\,-\, \frac{ w'(t)} {\int_\R {w'}^2}\int_\R w' (\zeta)\, \phi(y,\zeta)\, \mathrm{d} \zeta,
$$
also satisfies $L({\bar\phi} ) = 0$ and, in addition,
\be
\int_\R   w' (t)\, {\bar\phi}(y,t)\, \mathrm{d} t  = 0\foral y\in \R^{N-1} .
\label{orti}\ee
Now, the function
$$
\vp (y) := \int_\R \big|{\bar\phi}(y,t)\big|^2\, \mathrm{d}t,
$$
is well defined and smooth. We compute
$$
\Delta_y \vp (y) \,=\, 2 \int_\R \Delta_y {\bar\phi} \cdot {\bar\phi} \,\mathrm{d}t
\,+\, 2\int_\R \big|\nabla_y {\bar\phi} \big|^2 \,\mathrm{d}t,
$$
and hence
\begin{align}
\begin{aligned}
0 &= \int_\R  \big(L({\bar\phi}) \cdot {\bar\phi} \big)
\\
& =
\frac{1}{2}\Delta_y \vp
- \int_\R  \big|\nabla_y {\bar\phi} \big|^2\, \mathrm{d}t
-\int_\R \big(\, |{\bar\phi} _t|^2 - F'(w){\bar\phi} ^2\,\big)\, \mathrm{d}t\, .
\end{aligned}
\label{ew}
\end{align}

From \equ{orti} and \equ{q}, we then get
$\frac{1}{2}\Delta_y \vp   -\gamma \vp   \ge 0. $
Since $\vp$ is bounded, it must be
 zero. In particular this implies that
 the bounded function
$$
g(y) = \int_\R w_\zeta (\zeta)\, \phi(y,\zeta)\, \mathrm{d} \zeta
$$
is harmonic and bounded, hence a constant. We conclude that
 $\phi(y,t) = Cw'(t)$, as desired.\qed

\bigskip

\noindent{\bf Proof of Proposition \ref{prop1}: }
We begin by proving  a priori estimates.
Let   $0<\sigma< \sqrt{2}$.
We first claim that there exists a constant $C>0$ such that for all small $\ve $
and every solution $\phi$ to Problem $\equ{p1}$ with
 $\|\phi\|_{\infty, \nu, \sigma} <+\infty$
and right hand side $g$ satisfying
 $\|g\|_{p,\sigma} < +\infty $ we have
\be
\| D^2 \phi\|_{p,\sigma } + \|D\phi\|_{\infty , \sigma} + \|\phi\|_{\infty , \sigma } \le  C\|g\|_{p , ,\sigma}.
\label{1.1}\ee
To establish this fact, it clearly suffices to consider 
the case $c(y)\equiv 0$.
By local elliptic estimates, it is enough to show that
\be
\|\phi\|_{\infty,\sigma } \le   C\|g\|_{p , \sigma}.
\label{2}\ee
Let us assume by contradiction that \equ{2} does not hold. Then we have sequences $\ve =\ve _n\to 0$,
$g_n$ 
with  $\|g_n\|_{p,\sigma}\to 0$, $\phi_n$ with $\|\phi_n\|_{\infty , \sigma} =1$  such that

\begin{align}
\begin{aligned}
 \partial_{tt}\phi_n + \Delta_{ \KK_\ve } \phi_n   +F'(w(t))\phi_n  =  g_n     \quad\hbox{in }   \KK_\ve \times \R,
\\
\int_\R \phi_n(y,t)\,w'(t)\,\mathrm{d}t=0  \foral  y\in  \KK_\ve  \, .
\end{aligned}
\label{p12}
\end{align}
Then we can find points $(p_n,t_n)\in  \KK_\ve \times \R$ such that

$$
 e^{-\sigma|t_n|}\,  |\phi_n(p_n,t_n)| \ge \frac 12.
 $$

We will consider different possibilities.  
Let us consider the local coordinates for $\KK_{\ve _n}$ around $p_n$,

$$
Y_{p_n, \ve _n} (\py)= \ve _n^{-1} Y_{\ve _n p_n} (\ve _n \py), \quad |\py|< \frac 1 {\ve _n},
$$
where $Y_p(\py)$ is given by \equ{Yp}.
Let us assume first that
$ |t_n| \le C .$
Then, the Laplace-Beltrami operator of $\KK_{\ve _n}$ takes locally the form
$$
\Delta_{\KK_{\ve _n}} \ = \ a_{ij}^0 (\ve _n \py )\pp_{ij}
 + \ve _nb_j^0 (\ve _n \py )\pp_j
$$
where
$$
a_{ij}^0 (\ve _n \py ) = \delta_{ij} + o(1),  \quad     b_i^0 (\ve_n \py) = O(1).
$$

Thus

$$
a_{ij}^0 \pp_{ij}\ttt \phi_n
 + \ve _nb_j^0\pp_j\ttt \phi_n   +  \partial_{tt}\ttt \phi_n  + F'(w (t) )\ttt \phi_n =  \ttt g_n (\py ,t),
 \quad |\py|< \frac 1\ve,
$$
where   $\ttt g_n (\py ,t): =  g_n\big(  Y_n( \ve  \py) , t\big)$.
We observe that this expression is valid for $\py$ inside the domain $\ve ^{-1}{\mathcal U}_k$ which is expanding to entire $\R^{N-1}$.
Since $\ttt\phi_n$ is bounded, and $\ttt g_n \to 0$ in $L^p_{loc}(\R^N)$, we obtain local uniform $W^{2,p}$-bounds.  Hence we may assume, passing to a subsequence, that
$\ttt \phi_n $ converges uniformly in the compact subsets of $\R^N$ to a function $\ttt \phi (\by, t)$ that satisfies

$$
\Delta_{\R^{N-1}}\ttt \phi
   +  \partial_{tt}\ttt \phi  + F'(w (t) )\ttt \phi
= 0 \, .
$$
Thus $\ttt \phi$ is non-zero and bounded.  Hence Lemma \ref{lemma l1}
implies that, necessarily,
$\ttt \phi (\by,t) = Cw'(t)$. On the other hand, we have
$$
0= \int_\R \ttt \phi_n (\by,t) \, w'(t)\, \mathrm{d}t \longrightarrow \int_\R \ttt \phi (\by,t) \, w'(t)\, \mathrm{d}t\
\quad\hbox{as } {n\to \infty}.
$$
Hence, necessarily $\ttt\phi \equiv 0$. But  $|\ttt \phi_n ( 0, t_n)| \ge \frac 12 $, and  $t_n$ was bounded,
the local uniform convergence implies $\ttt \phi \ne 0$. We have reached a contradiction.

\medskip
Now, if  $t_n$ is
unbounded, say, $t_n\to +\infty$, the situation is similar. The
variation is that we define now
$$
 \ttt \phi_n (\by, t) =    e^{\sigma(t_n + t) } \phi_n ( \by, t_n+  t), \quad  \ttt g_n (\by, t) =    e^{\sigma(t_n + t) } g_n (  \by, t_n+  t).
 $$
 Then $ \ttt \phi_n $ is uniformly bounded, and $\ttt g_n \to 0$ in $L^p_{loc} (\R^N)$.
 Now $\ttt \phi_n$ satisfies 
$$
{a_{ij}^{0}} ( \ve _n \by)\, \pp_{ij} \ttt \phi_n
\, + \,\pp_{tt} \ttt \phi_n
\,+ \, \ve _n b_j ( \py_n + \ve _n \by)\, \pp_{j} \ttt \phi_n \
$$
$$
 - 2\sigma \,\pp_t \ttt \phi_n  \, + \,F'\big(w(t+t_n ) +\sigma^2 \big)\, \ttt \phi_n \,=\, \ttt g_n .
 $$
 We fall into the limiting situation
\be
\, \Delta_{\R^{N-1}}  \ttt \phi \, + \,\pp_{tt}^2 \ttt \phi  \,
 - \,2\sigma \,\pp_t \ttt \phi  \, - \,( 2 -\sigma^2 )\, \ttt \phi \,=\,  0  \quad \hbox{in }\R^N
\label{positive}\ee
with  $\ttt \phi\ne 0$ bounded.
The  maximum principle implies that $\ttt \phi \equiv 0$. We obtain a contradiction that proves the validity of \equ{1.1}.

\medskip
It remains to prove existence of a solution $\phi$ of problem \equ{p1} with $\|\phi\|_{\infty,\sigma}< +\infty$.
We assume first that $g$ has compact support.
For such a $g$, Problem \equ{p1} has a variational formulation.
Let
$$
{\mathcal H} = \Bigg\{ \phi \in H^1_0( \KK_\ve \times \R) \ /\  \int_\R \phi(y,t)\,w'(t)\,\mathrm{d}t =0  \foral  y\in  \KK_\ve \, \Bigg\} \, .
$$
${\mathcal H}$ is a closed subspace of $H_0^1( \KK_\ve \times \R) $, hence a Hilbert space when endowed with  its natural norm,
$$
\|\phi\|_{\mathcal H}^2 = \int_{\KK_\ve} \int_\R \Big(\, |\pp_t\phi|^2 + |\nn_{\KK_\ve } \phi |^2 - F'(w(t))\, \phi^2 \, \Big)\,
\mathrm{d}V_{\KK_\ve} \, \mathrm{d}t\ .
$$
$\phi$
is then a weak solution of Problem \equ{p1} if
$\phi\in {\mathcal H}$ and satisfies
\begin{align*}
a(\phi, \psi) &:= \int_{\KK_\ve \times \R} \Big(\, \nn_{\KK_\ve } \phi \cdot \nn_{\KK_\ve } \psi
\, -\, F'(w(t))\, \phi\, \psi\,\Big)\, \mathrm{d}V_{\KK_\ve} \, \mathrm{d}t\,
\\
&\,=
- \int_{\KK_\ve  \times \R}  g\, \psi\, \mathrm{d}V_{\KK_\ve} \, \mathrm{d}t \foral \psi\in {\mathcal H}.
\end{align*}
It is standard to check that a weak solution of Problem \equ{p1}  is also classical provided that $g$ is regular enough.
Let us observe that because of the orthogonality  condition defining ${\mathcal H}$ we have that
$$
\gamma \int_{\KK_\ve  \times \R} \psi^2 \, \mathrm{d}V_{\KK_\ve} \, \mathrm{d}t\ \le \ a(\psi, \psi) \foral \psi \in {\mathcal H}.
$$
Hence the bilinear form $a$ is coercive in ${\mathcal H}$,
and existence of a unique weak solution follows from Riesz's theorem.
If $g$ is regular and compactly supported,
$\psi$ is also regular. Local elliptic regularity implies in particular that $\phi$ is bounded.
Since for some $t_0>0$, the equation satisfied by $\phi$ is
\be
\Delta \phi + F'(w(t))\, \phi = c(y) w'(t), \quad |t|> t_0 , \quad y\in \KK_\ve,
\label{www}\ee
and $c(y)$ is bounded, then enlarging $t_0$ if necessary, we see that for $\sigma < \sqrt{2}$,
a suitable barrier argument shows that $|\phi| \le  Ce^{-\sigma|t|}$, hence
$\|\phi\|_{p,\sigma} < +\infty  $.  From \equ{1.1} we obtain that
\be
 \|D^2 \phi\|_{p,\sigma} + \|D \phi\|_{\infty  , \sigma} + \| \phi\|_{\infty  , \sigma} \le C \|g\|_{p , \sigma}.
\label{coota}\ee
Now,  for an arbitrary $\|g\|_{p,\sigma} <+\infty$ we consider a sequence of compactly supported  approximations uniformly controlled in $\|\ \|_{p , \sigma}$ (thus inheriting corresponding control on the approximate solutions).  Passing to a subsequence if necessary, we obtain local convergence to
a solution $\ttt\phi$ to the full problem which respects the estimate \equ{cota}. This concludes the proof. \qed

\setcounter{equation}{0}
\section{Solving the nonlinear projected problem}

To solve  Problem  \equ{glufinp} and for the subsequent step of adjusting $\hh$ so that the quantities
$c_l(y)$ are all identically zero, it is important to keep track of the Lipschitz character of the operators involved in this equation. We have the following result.

\begin{lemma}\label{lipNj}
There is a constant $C>0$ such that for all $\hh^l$ satisfying \equ{ass} and all $\phi^l$ with $\|\phi^l\|_{2,p,\sigma} \le \ve^{2-\tau}$, $l=1,2$
 we have
\be
\| N_j(\phi^1,\hh^1) -N_j(\phi^2,\hh^2)\|_{p,\sigma}
\, \le \,
C\ve^{4-\tau} \| \hh_1- \hh_2\|_{W^{2,p}(\KK)} +  \ve^{2-\tau} \| \phi^1- \phi^2\|_{2,p,\sigma},
\label{lipsNj}
\ee

\be
\| S_j(\hh^1) -S_j(\hh^2)\|_{p,\sigma} \ \le \ C\ve^{3-\tau} \| \hh_1- \hh_2\|_{W^{2,p}(\KK)}.
\label{lipSj}
\ee
\end{lemma}

\proof

We have to check the Lipschitz character of the operators $\NNN_j(\phi,\hh)$ in \equ{Nj} in the norm $\|\ \|_{p,\sigma}$. Let us consider
each of the terms in  formula \equ{Nj}.
Let us consider first the operator $\BB_j\phi_j$ in \equ{BBj}. On $\phi$ and $\hh$ we assume
\be
\|\phi\|_{2,p,\sigma} \le \ve^{2-\tau}, \quad \quad \|\hh\|_{W^{2,p}(\KK)} \le M.
\label{ass2}\ee
 This operator has the form in local coordinates
$$
\BB_j\phi_j (y,t) = \BB\big( h_j, \pp_{ik} h_j, \pp_i h_j,\phi_j, \py,t \big)  .
$$
Let us consider the operator $B_j\phi_j$ in \equ{Bj}.
We see that the explicit dependence on $h_j$ comes only from the coefficients
$a_{ik}$ and $b_i$, more precisely on smooth functions of the form
$a\big(\ve\py, \ve t + \ve f_j(\ve \py)\big)$, $f_j = \xi_j  + h_j$, so that
$\pp_{h_j} a =  O(\ve).
$
We also find
$$
\pp_{h_j}\chi_{j3} = O(1),
\quad
\pp_{h_j} F'(\ww_j) = \sum_{k\ne j} w'\big(t- (f_k -f_j)\big) = O(\ve^{2-\tau}).
$$
Taking these facts into account
we then find that for arbitrarily small $\tau >0$,
$$
\pp_{h_j} \BB_j\phi_j \,=\,  O(\ve^{1-\tau})  D_{\MM_\ve}^2\phi_j
\,+\,  O(\ve^{2-\tau}) D_{\MM_\ve}\phi_j
\,+\,O(\ve^{2-\tau}) \phi_j,
$$
and hence
$$
\|\pp_{h_j} \BB_j\phi_j\|_{0,p,\sigma} \le  C\ve^{1-\tau}\|\phi\|_{2,p,\sigma}\, .
$$
Observe that we have as well that
$$
\|\BB_j\phi_j\|_{0,p,\sigma} \le  C\ve^{1-\tau}\|\phi\|_{2,p,\sigma}\, .
$$
Let us consider the dependence on the derivatives of $h$.
We easily check that
$$
\pp_{D_\KK h} \BB_j\phi_j   =  O( \ve ) D^2\phi, \quad \pp_{D^2_\KK h} \BB_j\phi_j   = O(\ve ) D\phi_j.
$$
As a conclusion we find that, emphasizing the dependence on $\hh$ of the operator $\BB_j$,
\be
\|\BB_j(\phi^1 ,\hh^1) - \BB_j(\phi^2 ,\hh^2)\|_{0,p}  \le \ve^{1-\tau} \|\phi^1 -\phi^2\|_{2,p,\sigma} + \ve^{3-\tau}\, \|\hh^1 -\hh^2\|_{W^{2,p}(\KK)}.
\label{ppp}\ee
Let us consider the remaining operator in $\NNN_j$,
$$
\NN(\phi,\hh) :=
\chi_{j1}\, \left[  \,(F'(\ww_j) +2) \Psi_j(\phi,\hh) +  N(\Psi_j(\phi,\hh) + \phi_j) \,\right ].
$$
We write it as
$$
 \NN(\phi,\hh)(y,t) = \ttt \NN ( \phi, \psi, \hh,  y,t), \quad \psi = \Psi_j(\phi,\hh),
$$
and recall from Lemma \ref{lemapsi} that  $\|\psi\|_\infty = O(\ve^{4-\tau})$.
Observe first that
$$
\pp_\psi \ttt \NN  =  \chi_{j1}\, \left[  \,(F'(\ww_j) +2)  +  N'( \psi  + \phi_j) \,\right ] = O(\ve e^{-\sigma|t|} ),\quad
$$

$$
\pp_{\phi_j} \ttt \NN  =  \chi_{j1}\, N'(\psi + \phi_j)   = O( |\psi| + |\phi_j|) = O(\ve^{2-\tau}).
$$
In addition, we also check that
$$
\pp_{\hh} \ttt N =  O( |\psi|e^{-\sigma|t|}  + |\phi|^2)  = O( \ve^{4- \tau}).
$$
Using these estimates, and writing $\psi^l= \Psi_j(\phi^l,\hh^l)$ we find
$$
\|\NN(\phi^1,\hh^1) - \NN(\phi^2,\hh^2)\|_{p,\sigma}\  = \
\| \ttt N (\phi^1,\psi^l,\hh^1,\cdot) - \ttt N (\phi^2,\psi^2,\hh^2,\cdot )\|_{p,\sigma}
$$
$$
\,\le\,C\ve\|\psi^1 -\psi^2\|_\infty +  C \ve^{2-\tau}  \|\phi^1 -\phi^2\|_{2,p,\sigma} +  C\ve^{4-\tau}\|\hh^1 -\hh^2\|_{\infty}.
$$
Recalling now, \equ{lippsi}, and combining this with estimate \equ{ppp} we arrive to the desired result. The proof of \equ{lipNj} is concluded.
The proof of estimate \equ{lipSj} is similar, taking into account the explicit form of the error.
\qed

\medskip
\begin{prop} \label{pop} Given $\hh$ satisfying \equ{assh},
problem \equ{glufinp} has a unique solution $\phi = \Phi(\hh)$ with $\|\phi\|_{2,p,\sigma} \le \ve^{2-\tau}$. Moreover, we
have the validity of the Lipschitz conditions
\be
\|\Phi(\hh^1) - \Phi(\hh^2)\|_{2,p,\sigma} \ \le \ C\,\ve^{2-\tau} \, \|\hh^1- \hh^2\|_{W^{2,p}(\KK)} .
\label{i1}\ee
In addition, we have that
\be
\|\NNN_j(\Phi(\hh^1),\hh^1) - \NNN_j(\Phi(\hh^2),\hh^2)\|_{2,p,\sigma} \le  C\,\ve^{4-\tau} \, \|\hh^1- \hh^2\|_{W^{2,p}(\KK)}.
\label{i2}\ee

\end{prop}

\proof
Let $T(g)$ be the operator defined as the solution of \equ{p1} predicted by Proposition \ref{prop1}. Then we find a solution to \equ{glufinp}
if we solve the fixed point problem for $\phi=(\phi_1,\ldots, \phi_N)$
\be
\phi_j = \mathfrak{B}_j(\phi,\hh):=
- T\big(S_j (\hh ) + \NNN_j(\phi,\hh) \big) \foral j=1,\ldots, N.
\label{fp2}\ee
We will check that the operator
$\mathfrak{B}(\phi,\hh) = \big(\mathfrak{B}_1(\phi,\hh),\ldots, \mathfrak{B}_N(\phi,\hh)\big)$
is a contraction mapping in
$\phi$ in a ball for the norm $\|\ \|_{2,p,\sigma}$. We will do more, checking as well the Lipschitz dependence  in $\hh$.
Using the above lemma we conclude that the operator $\mathfrak{B}$ is a contraction mapping on the region $\|\phi\|_{2,p,\sigma}\le \ve^{2-\tau}$.
Now, using \equ{psi1},
$$
|\mathfrak{B}(0)|
\,\le\,
\chi_{j1}\Big(\, \big(F'(\ww_j) + 2\big)|\Psi(0)|   +   |\Psi(0)|^2 +   |S_j (\ww )|\Big) .
$$
Thus
$$
|\mathfrak{B}(0)|
\,\le\,
\ve^{4-\tau} e^{-\sigma|t|} +  \chi_{j1}\ve^{7} + C\ve^2|D^2_{\KK}h(\ve y)|e^{-\sigma|t|} + C\ve^2e^{-\sigma|t|},
$$
and hence
$$
\|\mathfrak{B}(0)\|_{p,\sigma}\le  C\ve^{2-\tau} .
$$
As a conclusion, we can apply the contraction mapping principle, and find a unique solution $\phi$ of problem \equ{fp2} such that
$$
\| \phi\|_{2,p,\sigma} \le \beta \ve^{2-\tau},
$$
for a suitably large choice of $\beta$.

\setcounter{equation}{0}
\section{The Jacobi-Toda system}

Once problem \equ{glufinp} has been solved by $\phi = \Phi(\hh)$, according to Proposition \ref{pop},
the remaining task is to find an $\hh$ such that for all $\ell =1,\ldots, m$, we have
\be
 I_\ell (y) = {\int_\R \big( S_\ell (\hh ) + \NNN_j(\Phi(\hh),\hh)\big)\, w'\,\mathrm{d}t }\, = 0 \foral y\in \KK_\ve.
\label{equc}\ee
%
Using the definition of $S_\ell$ in \equ{Sjj}, expansion \equ{Sn}, Lemma \ref{error1} and the definition of $\NNN_j(\Phi(\hh),\hh)$, we get
\begin{align}
\ve^{-2}I_{\ell}( \ve^{-1} y) \,=\ &b_1\,  \Bigl(\Delta_{\KK} h_{\ell} +   K( y) f_{\ell} ,\Bigr)
\,-\,  b_2\rho_\ve\,\left[\,  e^{ -\sqrt{2}(h_{\ell } -h_{\ell -1})}
-  e^{ -\sqrt{2}(h_{\ell +1} -h_{\ell }) } \, \right ]\nonumber
\\
&\,+\, \theta_\ell(\hh),
\label{proyec}
\end{align}
where  $\theta_\ell $ is a small operator:

$$
\| \theta_\ell (\hh)\|_{L^p(\KK)}\, = \, O(\ve^{1-\tau} ),
$$
 for any $\tau > \frac {N-1}p$, uniformly on $\hh$.
The constants $b_1,b_2$ are given by
$$
b_1 = \int_\R {w'(t)^2}\, \mathrm{d}t,
\quad
b_2  = \int_\R  6(1- w^2(t))e^{ \sqrt{2} t} w'(t) \mathrm{d}t
= \int_\R  6(1- w^2(t))e^{ -\sqrt{2} t} w'(t) \mathrm{d}t.
$$
\\

\noindent{\bf Part I: }
Recall the relation in (\ref{fj})
$$
f_{\ell}(y) = \Big(\ell -\frac{m+1}2\Big)\rho_\ve + h_{\ell}(y).
$$
Since we want that the functions $h_{\ell}$  make the quantities
 $I_{\ell}$ as small as possible, it is reasonable to find first an $\hh$ such  that the equations,
for $\ell =1,\ldots m$,
\be
 b_1\,  \Bigl(\Delta_{\KK} h_{\ell} +   K(y) f_{\ell} ,\Bigr)
 \,  -\,  b_2 \rho_\ve\,\left[\,  e^{ -\sqrt{2}(h_{\ell } -h_{\ell -1})} -   e^{ -\sqrt{2}(h_{\ell +1} -h_{\ell }) } \, \right ] = 0,
\label{pjj}
\ee
be approximately satisfied. We set
\be
R_{\ell}(\hh) :=  \sigma  \Bigl(\Delta_{\KK} h_{\ell} +   K(y) f_{\ell} ,\Bigr)\,  -\,  \,\left[\,  e^{ -\sqrt{2}(h_{\ell } -h_{\ell -1})} -   e^{ -\sqrt{2}(h_{\ell +1} -h_{\ell }) } \, \right ] ,
\ee
where
$$\sigma:= \sigma_\ve = \rho_\ve^{-1} b_1 b_2^{-1} \sim (\log \frac 1\ve )^{-1}\ ,$$
and
\be
\bR(\hh) :=  \left [\begin{matrix} R_1(\hh)\\ \vdots \\ R_m(\hh)  \end{matrix}   \right ] .
\label{R}
\ee

We would like  to find a solution $\hh $ to the system $ \bR(\hh )=0$.
To this end,
we find first a convenient representation of the operator $\bR(\hh)$. Let us consider the auxiliary variables
$$ \vv :=  \left [\begin{matrix} \bar \vv\\ v_m  \end{matrix}   \right ], \quad \bar\vv = \left [\begin{matrix} v_1\\ \vdots \\ v_{m-1}  \end{matrix}   \right ],$$
defined in terms of $\hh$ as

$$
v_{\ell} = h_{\ell +1} -h_{\ell}\quad\mbox{with } \ell=1,\ldots, m-1,
\quad
v_m = \sum_{\ell =1}^m h_{\ell},
$$
with the conventions $v_0= v_{m+1} =+\infty$ and define  the operators
$$
\bS(\vv) :=  \left [\begin{matrix} \bar \bS(\bar \vv) \\ S_m( v_m)  \end{matrix}   \right ],
\quad
\bar\bS(\bar\vv) = \left [\begin{matrix} S_1(\bar\vv)\\ \vdots \\ S_{m-1}(\bar\vv) \end{matrix}   \right ].
$$
where we have setted
$$
S_{\ell}({\tt v} ) : = R_{\ell +1}(\hh) - R_{\ell}(\hh)\ =\
$$

$$
\sigma \, \Bigl(\Delta_{\KK} v_{\ell} +   K( y) (\rho_\ve + v_{\ell})\Bigr)\,  + \, \left\{ \begin{matrix}
  e^{ -\sqrt{2}v_{2}} -   2e^{ -\sqrt{2}v_{1} }  \,  & \hbox{if}& \ell= 1,
\\
      e^{ -\sqrt{2}v_{\ell +1}} -   2e^{ -\sqrt{2}v_{\ell} } +  e^{ -\sqrt{2}v_{\ell -1}} \,   & \hbox{if}& 1<\ell<m-1,
\\
  -   2e^{ -\sqrt{2}v_{m-1} } +  e^{ -\sqrt{2}v_{m -2}} \,   & \hbox{if}& \ell=m-1, \\ \end{matrix}\right.
$$

\noindent
and

$$
S_m( v_m ) := \sum_{\ell =1}^m R_{\ell}(\hh) =   \sigma\, \Bigl(\Delta_{\KK} v_m +   K( y) v_m  \Bigr).
$$
Then the operators $\bR$ and $\bS$ are in correspondence through the formula
\be
\bS({\tt v}) \ = \ {\bf B}\, \bR\left(\,{\bf B}^{-1} {\tt v} \, \right ),
\label{bS}\ee
where ${\bf B}$ is the constant, invertible $N\times N$  matrix
\begin{align}\label{definitionofbfB}
 {\bf B}\, =\,
\left[\begin{matrix}
-1 & \ 1& 0 & \cdots  &  0
\\
                      0 & -1&1 & \cdots   &  0
\\
                      \vdots   &  \ddots & \ddots & \ddots   & \vdots
\\
                       0 &    \cdots   & 0&  -1 &  1
\\
                       1 &    \cdots &1 &\ 1&  1
\\
\end{matrix}\right ],
\end{align}

\noindent and then the system $\bR(\hh)=0$ is equivalent to $\bS(\vv)=0$,
which setting $\beta = b_2b_1^{-1}$ decouples into

\begin{align}
\bar \bS(\bar \vv ) \, =\,
\sigma \Bigl[ \Delta_{\KK} {\bar \vv}  +   K( y) {\bar \vv} \Bigr ] \, + \, \beta K( y)\, \left[\begin{matrix} 1
\\
\vdots \\ 1 \end{matrix}\right ]\ + \ \bar\bS_0(\bar \vv ) \ = \ 0, \label{skk0}
\\
\nonumber
\\  S_m(v_m)\, =\,  \sigma \left (\Delta_{\KK} {v_m}  +   K( y) {v_m}\right ) \ =\ 0,\qquad\qquad
\label{skk}\end{align}
where
\be
  \bar \bS_0  (\bar \vv) :=   -{\bf C}\, \left[\begin{matrix} e^{-\sqrt{2}v_1} \\  \vdots \\ e^{-\sqrt{2}v_{m-1}}\end{matrix}\right ],\quad
{\bf C}\, =\,
\left[\begin{matrix} 2 &- 1& 0 & \cdots  &  0 \\
                      -1 & 2&-1 & \cdots   &  0 \\
                      \vdots   &  \ddots & \ddots & \ddots   & \vdots   \\
                       0 &    \cdots   &-1  & 2 & - 1 \\
                       0 &    \cdots & & -1&  2 \\
                          \end{matrix}\right ].
\label{bS0}\ee
In system \equ{skk0}-\equ{skk}, the second relation and our non-degeneracy assumption force $v_m=0$.
Thus we look for a solution $\vv =(\bar\vv ,0)$ of the system, where $\bar \vv$ satisfies \equ{skk0}.
Rather than finding an exact solution $\bar\vv$ of $\bar\bS(\bar \vv)=0$ we will find a good approximation. More precisely,
by means of a simple iterative procedure, we will find for each $k\ge 1$ a function $\bar\vv^k$ with the property that
\be
  \bar\bS(\bar \vv^k)= O(\sigma^{k} ). \label{conc}\ee
  Let us  find a function $\bar \vv^1$ with the desired property \equ{conc} for $k=1$.
We consider  the vector $\bar \vv^1(y)$ defined by the relations
$$
\bar \bS_0 (\bar \vv^1) \ =\ -{\bf C}\, \left[\begin{matrix} e^{-\sqrt{2}v_1^1} \\  \vdots \\ e^{-\sqrt{2}v^1_{m-1}}\end{matrix}\right ] =   -\beta K( y)\, \left[\begin{matrix} 1\\   \vdots \\ 1\end{matrix}\right ]\ .
$$
We  compute explicitly
\be
v_{\ell}^1(y)  =    \frac 1{\sqrt{2}}  \log\, \left [ \, \frac \beta 2 \,K(y)\, ( m-\ell)\,\ell\, \right ] , \ell =1,\ldots, m-1,
\label{v0}\ee
and get from \equ{skk0}
$$
\bar \bS ( \bar \vv^1)\ =    \  \sigma\,\left [ \Delta_{\KK} \bar \vv^1 +   K( y) \bar \vv^1   \right] \ =\ O(\sigma).
$$

This approximation can be improved to any order in powers of $\sigma$, as the following lemma states.

\begin{lemma}\label{lemita1}

Given $k\ge 1$,
there exists a function of the form
$$
\bar\vv^k(y,\sigma) = \bar \vv^1(y) +  \sigma \xi_k (y,\sigma), \quad
$$
where $\bar \vv^1(y)$ is defined by \equ{v0},  $\xi_1\equiv 0$,
and $\xi_k$ is smooth on $ \KK\times [0,\infty)$,
such that

$$\bar\bS( \bar\vv^k) = O(\sigma^k)$$
as $\sigma \to 0$, uniformly on $\KK$. In particular, $$\hh^k:= {\bf B}^{-1}\left [\begin{matrix} \bar\vv^k\\ 0\end{matrix}\right ] ,$$
 with ${\bf B}$ is given by
\equ{definitionofbfB}, solves approximately system \equ{pjj} in the sense that
 $$\bR( \hh^k) = O(\sigma^k). $$
\end{lemma}

\proof

In order to find a subsequent improvement of approximation beyond $\vv^1$,
we set $\bar\vv^2 = \bar\vv^1 + \omega_1$.
Let us expand
\be
\bar \bS(\bar \vv^1+ \omega ) \,  =\,
\sigma \Bigl[ \Delta_{\KK} {\omega}  +   K( y) {\omega} \Bigr ]
\, + \,\sigma\big(\Delta_{\KK} {\tt v}^1  +   K( y) {\tt v}^1\big)
\,+\, D \bar \bS_0(\bar \vv^1 )\omega
\,+\, \bN(\omega) ,   \quad
\label{sk0}
\ee
where
$$
 D\bar \bS_0( \bar \vv^1) \,  = \, \sqrt{2} {\bf C}\, \left [\begin{matrix}   e^{-\sqrt{2} v_1^1} &  0& \cdots & 0\\
                                                                                0 &   e^{-\sqrt{2} v_2^1}& \cdots & 0\\
                                                                                 \vdots &                   & \ddots &  \vdots  \\
                                                                                0   &  0& \cdots & e^{-\sqrt{2} v_{m-1}^1}
                                                                                \end{matrix} \right ]
\qquad\quad\qquad\qquad\qquad
$$

\medskip
\medskip

\be
\qquad =\,
\frac \beta  {\sqrt{2}}  \, K(y)\,
\left[\begin{matrix} 2a_1 & -a_2&  0 & &\cdots  &  0
\\
-a_1 & 2a_2& -a_3 & &\cdots   &  0
\\
0    &   -a_2 &  2a_3 & &\cdots & 0
\\
\vdots   &  &\ddots & \ddots & \ddots   & \vdots
\\
0 &   & \cdots   & -a_{m-3} & 2a_{m-2} &  -a_{m-1}
\\
0 &    &\cdots & & -a_{m-2}&  2a_{m-1}
\\
\end{matrix}\right ],
\label{E}
\ee
with
\be
a_{\ell} = \, ( m-\ell )\,\ell,\, \quad \ell =1,\ldots, m-1, \label{ann}
\ee
and
$$ \bN        (\bar\vv )\, = \, \frac \beta  {{2}}\, {\bf C}\, \left [\begin{matrix}   a_2(\,e^{-\sqrt{2} v_1^0} - 1 +\sqrt{2}\, v_1^0 \,)  \\
                                                                                 \vdots \\
                                                                                 a_m(\,e^{-\sqrt{2} v_m^0} - 1+  \sqrt{2}\, v_m^0 \,)  \\
                                                                                \end{matrix} \right ] \, .
$$

The matrix $D\bar \bS_0( \bar \vv^1)  $ is clearly invertible.
Let us consider the unique solution $\omega_1 = O( \sigma)$ of the linear system

\be
D\bar \bS_0( \bar \vv^1) \omega_1
\ =\  - \sigma\big(\Delta_{\KK} \bar{\tt v}^1  +   K( y) \bar {\tt v}^1\big)
\ =\ O(\sigma),
\label{om1}
\ee
and define $\bar \vv^2 = \bar\vv^1 +\omega_1$.
Then from \equ{sk0} we have

\be
\bar\bS(\bar\vv^2  )
\, =\,\sigma\big(\Delta_{\KK} \omega_1  +   K( y) \omega_1\big)   +   \bN( \omega_1)
\, =\, O( \sigma^{2})  . \quad
\label{sk1-1}\ee

\noindent Next we define $\bar\vv^3 = \bar \vv^2 + \omega_2$ where  $\omega_2 =  O( \sigma^{2})$
is the unique solution of

\be
 - D\bar \bS_0( \bar \vv^1) \omega_2 \,=\, \sigma\big(\Delta_{\KK} \omega_1  +   K( y) \omega_1\big)
  +   \bN( \omega_1).
\label{om2}\ee
Then from \equ{sk0} we get
\be
\bar\bS(\bar\vv^3  )
\,  =\,\sigma\big(\Delta_{\KK} \bar \omega_2  +   K( y)\bar  \omega_2\big)   +   \bN( \omega_1 + \omega_2)-\bN( \omega_1)
\, =\, O( \sigma^{2})  . \quad
\label{sk1}\ee

\medskip
In general, we define inductively, for $k\ge 3$,
$\bar \vv^{k+1} =\bar\vv^k +\omega_k$ where $\omega_k$ is the unique solution of
the linear system

\begin{align}
-D\bar \bS_0(\bar\vv^1)\, \omega_{k}
\,=\ &  \sigma\big(\Delta_{\KK} {\omega}_{k-1}  +   K( y) {\omega}_{k-1}\big)
\,+\,\bN ( \omega_1+\cdots + \omega_{k-1})\nonumber
\\
\nonumber
\\
&\,-\, \bN( \omega_1+\cdots + \omega_{k-2}).
\label{omk}
\end{align}

\noindent Then clearly $\omega_k= O(\sigma^k)$.
Let us estimate the size of $\bar\bS(\bar\vv^{k+1})$.
From \equ{sk0} we have

$$
\bar \bS(\bar \vv^{k+1} )
\  = \  \sigma\big(\Delta_{\KK} \bar{\tt v}^1  +   K( y) \bar{\tt v}^1\big)
\,+\, \Big[ \sigma(\Delta_{\KK}  + K  ) +  D\bS_0(\bar \vv^1)
 + \bN    \, \Big] \, \big(\,  \sum_{i=1}^{k} \omega_i \,\big).
$$

Now, using \equ{om1}, \equ{om2} and \equ{omk} we get

\begin{align*}
 &\Big[ \sigma(\Delta_{\KK}  + K  ) +  D\bS_0(\bar \vv^1) \Big]\, \big(\, \sum_{i=1}^{k} \omega^i\,\big)
\\
&\quad\,=\
\sigma (\Delta_{\KK} \bar\vv^1 + K \bar  \vv^1 )\, +  \,D\bar\bS_0(\bar \vv^1)\omega^1
\,+\,\sigma(\Delta_{\KK} \omega_ {k} + K \omega_{k} )
\\
&\qquad\, + \,\sum_{i=2}^{k} \left [ \sigma(\Delta_{\KK} \omega_ {i-1} + K \omega_{i-1} ) +  D\bar\bS_0(\bar \vv^1)\omega _i \right ]
\\
&\quad\, =\,- \bN     ( \omega^1)  -  \sum_{i=3}^k [\bN     (\omega_  1 +\cdots +\omega_ {i-1})- \bN     (\omega_ 1 +\cdots +\omega_ {i-2})]
\\
&\quad\, =\,
 - \bN(\omega_ 1 +\cdots +\omega_ {k}) \, .
\end{align*}

\noindent Hence,

\begin{align}
\bar\bS( \bar\vv^{k+1} ) \,  =\ & \sigma (\Delta_{\KK} \omega_{k} + K \omega_k )
\,+\, \bN     ( \omega_1 +\cdots +\omega_{k-1}+\omega_{k})
\,-\,\bN     (\omega_1 +\cdots +\omega_{k-1})\nonumber
\\
\,=\ &  O( \sigma^{k+1}).
\label{sk2}
\end{align}

\noindent Finally, the functions $\xi_1\equiv 0$ and

$$ \xi_k := \sigma^{-1}( \omega_1 +\cdots + \omega_{k-1}) ,\quad k\ge 2,
$$

\noindent
clearly satisfy the conclusions of the lemma, and the proof is concluded. \qed
\\

\noindent{\bf Part II:}
The question now, is how to use the approximation $\hh^k$ just constructed to find
an  exact $\hh$ solution to system \equ{equc}.
This system takes the form
\be
\bR( \hh) =  g,
\label{equhh}\ee
where $g$ is  a small function, actually a small nonlinear operator in $\hh$.
For the moment we will think of $g$ as a fixed function.
Since the operator $\bR$ decouples as in \equ{bS} when expressed in terms of $\bS$,
it is more convenient to consider
the equivalent problem

 \be
\bS( \vv) =  g,
\label{equvv}
\ee
which, according to expressions \equ{skk0} and \equ{skk}, decouples as
\begin{align}
\bar \bS(\bar \vv ) \, =\,
\sigma \Bigl[ \Delta_{\KK} {\bar \vv}  +   K( y) {\bar \vv} \Bigr ] \, + \, \beta K( y)\, \left[\begin{matrix} 1\\   \vdots \\ 1 \end{matrix}\right ]\ + \
  \bar\bS_0(\bar \vv ) \ = \ \bar g , \label{skkk0}
\\
\nonumber
\\
\bS_m(v_m)\, =\,  \sigma \left (\Delta_{\KK} {v_m}  +   K( y) {v_m}\right ) \ =\ g_m.
\label{skkk}
\end{align}

Equation \equ{skkk} has a unique solution $v_m$ for any given function $g_m$,
thanks to the nondegeneracy assumption.
Therefore we will concentrate in solving Problem \equ{skkk0}, for a small given $\bar g$.
Let us write
$$
\bar \vv = \bar\vv^k +  \omega,
$$
where $\bar\vv^k$ is the approximation in Lemma \ref{lemita1}.
We express \equ{skkk0} in the form
\be
\ttt L_\sigma (\omega)\,:= \,  - \sigma \Bigl[ \Delta_{\KK} \omega  +   K( y) \omega  \Bigr ] - D\bar \bS_0( \bar\vv^k)\omega
\, =\, \bar \bS(  \bar\vv^k )  +  \bN_1(\omega ) +\bar g,
\label{equp}
\ee
where
\be
 \bN_1(\omega ) := \bar \bS_0(\bar\vv^k + \omega )- \bar \bS_0(\bar\vv^k)  - D\bar \bS_0( \bar\vv^k) \omega  \,,
\ee
and
$\bS_0$ is the operator in \equ{bS0}.

The desired solvability theory will be a consequence of a suitable invertibility statement for the linear operator $\ttt L_\sigma$.  Thus we consider
the equation
\be
\ttt L_\sigma(\omega) = \ttt g \quad\hbox{in }\KK.
\label{equp1}\ee

This operator is vector valued. It is convenient to express it in self-adjoint form by replacing the matrix
$D\bar \bS_0( \bar \vv^k)$ with a symmetric one. We recall that we have

$$
D\bar \bS_0( \bar \vv^k)  \,  = \, \sqrt{2} {\bf C}\, \left [\begin{matrix}   e^{-\sqrt{2} v_1^k} &  0& \cdots & 0\\
                                                                                0 &   e^{-\sqrt{2} v_2^k}& \cdots & 0\\
                                                                                 \vdots &                   & \ddots &  \vdots  \\
                                                                                0   &  0& \cdots & e^{-\sqrt{2} v_{m-1}^k}
                                                                                \end{matrix} \right ],
$$
where
the matrix ${\bf C}$ is given in \equ{bS0}. ${\bf C}$ is symmetric and positive definite. Indeed, a straightforward computation
yields that its eigenvalues are explicitly given by
$$ 1,\, \frac 12, \ldots, \frac {m-1}m. $$
We consider the symmetric, positive definite square root matrix of ${\bf C}$  and denote it by ${\bf C}^{\frac 12}$ .
Then setting $$\omega := {\bf C}^{\frac 12}\psi,\quad g:=  {\bf C}^{-\frac 12}\ttt g,$$
 we see that equation \equ{equp1} becomes
\be
L_\sigma (\psi) := - \sigma \Delta_\KK \psi - {\bf A}(y,\sigma)\psi\ =\  g   \quad\hbox{in }\KK,
\label{bvp}\ee
where ${\bf A}$ is the symmetric matrix
\be
{\bf A}(y,\sigma)\, =\,  \sigma K(y)\, {\bf I}_{m-1}\, + \, \sqrt{2}\, {\bf C}^{\frac 12}\,\left [\begin{matrix}   e^{-\sqrt{2} v_1^k} &  0& \cdots & 0\\
                                                                                0 &   e^{-\sqrt{2} v_2^k}& \cdots & 0\\
                                                                                 \vdots &                   & \ddots &  \vdots  \\
                                                                                0   &  0& \cdots & e^{-\sqrt{2} v_{m-1}^k}
                                                                                \end{matrix} \right ]\, {\bf C}^{\frac 12}.
\ee

Since
$$
\vv^k = \vv^1(y) + \sigma \xi^k(y,\sigma),
$$
we have that $A$ is smooth in its variables and
\be
{\bf A}(y,0)\, =\,   \, \frac\beta {\sqrt{2}}\,K(y)\, {\bf C}^{\frac 12}\,\left [\begin{matrix}   a_1 &  0& \cdots & 0\\
                                                                                0 &   a_2& \cdots & 0\\
                                                                                 \vdots &                   & \ddots &  \vdots  \\
                                                                                0   &  0& \cdots & a_{m-1}
                                                                                \end{matrix} \right ]\, {\bf C}^{\frac 12}.
\ee
where $a_\ell = \ell(m-\ell)$. In particular, ${\bf A}(y,\sigma)$ has uniformly positive eigenvalues whenever $\sigma$ is
sufficiently small.

\medskip

\medskip
Our main result concerning uniform solvability of Problem \equ{bvp} is the following.

\begin{prop} \label{solving}

There exists a sequence of
values $\sigma = \sigma_\ell \to 0$ such that
$L_{\sigma}$ is invertible. More precisely, for any $g\in L^2(\KK)$  there exists a unique solution $\psi= L_{\sigma}^{-1} g\in H^1(\KK)$ to  equation \equ{bvp}.
This solution satisfies
\be
\sigma\|D^2_{\KK}{\psi}\|_{L^2(\KK)}+  \| \psi\|_{L^2(\KK)} \ \le \ C \sigma^{-\frac {N-1} 2} \| g\|_{L^2(\KK)} .
\label{k 15}\ee
Moreover, if $p>{N-1} $, there exist $C,\nu >0$ such that the solution satisfies besides the estimate
\begin{align*}
\|D^2_{\KK}{\psi}\|_{L^p(\KK)}
+ \,\|
D_{\KK}{\psi}\|_{L^\infty(\KK)} +\|{\psi}\|_{L^\infty(\KK)} \,\leq\,
C\sigma^{-\frac{N-1}2 - \nu
}\,{\|{{g}}\|_{L^p(\KK)}}.
\end{align*}
In addition, for $N=2$, we have the existence of positive numbers $\nu_1,\nu_2,\ldots, \nu_{m-1}$ such that for all small $\sigma$ with
$$
|   \nu_i\sigma  -   j^2 | > c\sigma^{-\frac 12} \foral j\ge 1,\quad i=1,\ldots, m-1,
$$
for some $c>0$, then $\psi= L_{\sigma}^{-1} g$ exists and estimate \equ{k 15} holds.

\end{prop}

\medskip
We postpone the proof of this result in the last section.
Assuming its validity, we will use it to derive
a solvability statement for the Problem  \equ{equhh}, and to conclude the corresponding solvability of system \equ{equc},
hence that of Theorem \ref{main}.

\setcounter{equation}{0}
\section{Solving system \equ{equc}: Conclusion of the proof of Theorem \ref{main}}\label{section5}

\subsection{Solving Problem \equ{equhh}}
Here we refer to objects and notation introduced in the previous section.

\medskip
Because of the definition of $L_\sigma$, the statement of Proposition \ref{solving} holds as well for the operator $\ttt L_\sigma$ in
 equation \equ{equp1}. Choosing $\sigma$ as in the proposition, we write this equation as the fixed point problem
 \be
 \omega \, =\,  T(\omega) := \ttt L_\sigma^{-1} \, \left ( \bar \bS( \bar\vv^k) + \bar g  +  \bN_1(\omega)  \right ).
\label{fp22}\ee
By construction, we have that
$$
\|\bar \bS( \bar\vv^k)\|_{L^p(\KK)}  \le C\sigma^{k} .
$$
On the other hand, if $\|\omega\|_{L^\infty(\KK)}\le \delta $, with $\delta$  sufficiently small we also have that
$$
\| \bN_1(\omega) \|_{L^\infty(\KK)} \le C \delta^2,
$$
and in this region
$$
\| \bN_1(\omega_1 ) -  \bN_1(\omega_2 ) \|_{L^\infty(\KK)}\le C\delta  \|\omega_1 -\omega_2 \|_{L^\infty(\KK)}.
$$

We observe then that, for $\nu$ as in Proposition \ref{solving},
$$
\| T(\omega_1 ) -  T(\omega_1 ) \|_{W^{2,p}(\KK)}
\,\le\,
C\delta \sigma^{- \frac{N-1}2 - \nu}  \|\omega_1 -\omega_2 \|_{L^\infty(\KK)},
$$
and
$$ \| T(\omega )\|_{W^{2,p}}
\,\le\,
C \sigma^{- \frac{N-1}2 - \nu}\, ( \sigma^k + \delta^2 + \|g\|_{L^p(\KK)}\, ).
$$

Thus if we choose

$$
k> 2\big(\frac{N-1}2 +\nu\big),
$$
and  $g$ with
$$
\|g\|_{L^p(\KK)}  \le  \sigma^k,
$$
then the choice $\delta =  \mu \sigma^{ k- \frac{N-1}2 -\nu}$
with $\mu$ large and fixed yields, thanks  to the contraction mapping principle,   the existence of a unique solution
$\omega$  to Problem \equ{fp22}, with $$\|\omega\|_{W^{2,p}(\KK)} \le  \mu \sigma^{ k- \frac{N-1}2 -\nu}. $$
Let us call $\omega =: \Omega(g)$. Then, in addition $\Omega$ satisfies the Lipschitz condition
$$
\| \Omega (g_1)  -  \Omega(g_2 ) \|_{W^{2,p}(\Omega)} \, \le \, C\, \sigma^{- \frac{N-1}2 - \nu}\,\| g_1-g_2\|_{L^p(\KK)}
$$
for all $g_1, g_2$ with $\|g\|_{L^p(\KK)}  \le  \sigma^k$.
It follows that the equation \equ{equvv}
 \be
\bS( \vv) =  g
\label{equvv1}\ee
can be solved under these conditions.
In the form
$$ \vv = V(g):= \left[ \begin{matrix} \bar \vv^{k} + \Omega (\bar g ) \\ \sigma^{-1} ( \Delta_\KK + K )^{-1} g_m \end{matrix} \right ], $$
and  therefore  the equation
$
\bR( \hh) =   \ttt g$
can be solved for any given small $\ttt g\in L^p(\KK)$ by means of the correspondence $$\bS({\tt v}) \ = \ {\bf B}\, \bR\left(\,{\bf B}^{-1} {\tt v} \, \right ). $$
 This yields the following result.

\begin{lemma}\label{lemi3}
Given $k> 2( \frac{N-1}2 + \nu) $,
then for all sufficiently small $\sigma$ satisfying the statement of Proposition \ref{solving},
and all functions $\ttt g$ with

$$
\|\ttt g\|_{L^p(\Omega)} \le \sigma^{k} ,
$$
there exists a solution of the equation
\be
\bR( \hh) =   \ttt g,
\label{equhhh}
\ee
of the form
$$
\hh = \hh^k + H(\ttt g),
$$
where the operator $H$ satisfies
$$
\| H (\ttt g_1)\|_{W^{2,p}(\Omega)} \, \le \, C\, \sigma^{k - \frac{N-1}2 - \nu},
$$
and
$$
\| H (\ttt g_1) - H (\ttt g_2) \|_{W^{2,p}(\Omega)}
\, \le \,
C\, \sigma^{- \frac{N-1}2 - \nu}\,\| \ttt g_1-\ttt g_2\|_{L^p(\KK)}.
$$
\end{lemma}
\qed

\bigskip
\subsection{Proof of Theorem \ref{main}}
We need to prove the existence of $\hh$ satisfying System \equ{equc}.
According to expansion \equ{proyec},  we have that
\be
 \frac \sigma {b_2} \ve^{-2} I_{\ell}( \ve^{-1} y) \,=\, \bR(\hh ) \, -\, G(\hh),
 \ee
 where
 $$
 G(\hh) =
 -\sigma  \frac \sigma {b_2}\,\theta_\ell (y),
 $$
 and $\theta_\ell$ is the remainder in \equ{proyec}.
 We will estimate this operator. We have that
 $$
 \theta_\ell (y) =    \underbrace{(-1)^{\ell -1}\ve^{-2} \int_\R \zeta_\ve \Theta_\ell (\hh)( \ve^{-1}y,t) w'(t)\, \mathrm{d}t}_{Q_1(\hh)}
  + \underbrace{\ve^{-2} \int_\R \NNN_\ell(\Phi(\hh),\hh)( \ve^{-1}y,t)\, w'\,\mathrm{d}t}_{Q_2(\hh)} \,
 $$
 where $\Theta_\ell$ is described in \equ{Sn}.
We have, using  Lemma \ref{error1},
\begin{align*}
\| Q_1(\hh)\|_{L^p(\KK)}
\,\le\ &
\ve^{\frac{N-1}p-2}  \int_\R \| \zeta_\ve \Theta_\ell (\hh)( \cdot ,t)\|_{L^p(\KK_\ve)}\, w'(t)\, \mathrm{d}t
\\
\nonumber
\\
\,\le\ &
C \ve^{-2}\| \zeta_\ve \Theta_\ell (\hh)( \cdot ,t)\|_{p,\sigma}
\, \le\,  C\ve^{1- \tau} .
\end{align*}
And similarly, using  Lemma \ref{lispchitzerror} we get
$$
\| Q_1 (\hh^1)  -  Q_1(\hh^2) \|_{L^p(\KK)} \le  C\ve^{1-\tau} \|\hh^1-\hh^2\|_{W^{2,p}(\KK)}.
$$
whenever the vector-valued functions $\hh_1$, $\hh_2$ satisfy constraints \equ{assh}.
A similar argument, using  the estimates for the operator $N_j(\Phi(\hh),\hh)$ in Proposition \ref{pop} yields

$$
\| Q_2 (\hh) \|_{L^p(\KK)} \le  C\ve^{2-\tau} ,\quad
\| Q_2 (\hh^1)  -  Q_2(\hh^2) \|_{L^p(\KK)} \le  C\ve^{2-\tau} \|\hh^1-\hh^2\|_{W^{2,p}(\KK)}.
$$
As a consequence, the operator $G(\hh)$ satisfies

\be
\| G(\hh) \|_{L^p(\KK)} \le  C\ve^{2-\tau} ,\quad \|G(\hh^1) -G(\hh^2)\|_{L^p(\KK)}\,  \le \, C\ve^{1-\tau} \|\hh^1- \hh^2\|_{W^{2,p}(\KK)}, \quad
\label{lok}\ee

Thus we need to solve the system

\be
\bR ( \hh ) = G(\hh),
\ee
which can be rewritten in the form
\be
\bR ( \hh^k + {\tt q} ) =  G( \hh^k + {\tt q} ).
\label{fpN}\ee

\noindent We use the operator $H(g)$ defined in Lemma \ref{lemi3}, and look for a solution of \equ{fpN} by solving

\be
\hh =  H( G( \hh^k + {\tt q}  ) )=:  D( {\tt q} ),
\label{fpN1}\ee
for a sufficiently large $k$, in the region
$$
{\mathcal R} = \Big\{ {\tt q}  \in W^{2,p}(\KK)\ /\ \| {\tt q} \|_{W^{2,p}(\KK)} \le \mu\sigma^{k- \frac{N-1}2 -\nu}     \Big\},
$$
for a sufficiently large $\mu$.  From Lemma \ref{lemi3} and \equ{lok}, we get that

$$
\| D( {\tt q} _1 ) -  D( {\tt q} _2 ) \| _{W^{2,p}(\KK)} \, \le \, C\, \sigma^{- \frac{N-1}2 - \nu}\,\| G( \hh^k +  {\tt q} _1 )- G( \hh^k + {\tt q} _2 )\|_{L^p(\KK)}
$$
$$
\qquad\qquad\,\le\, C\, \sigma^{- \frac{N-1}2 - \nu}\,\ve^\tau \|{\tt q} _1 - {\tt q} _2\|_{W^{2,p}(\KK)},
$$
Hence $D$ is a contraction mapping in ${\mathcal R}$. Besides we have
$$
\| D(0)\| _{W^{2,p}(\Omega)} \, \le \, C\, \sigma^{k - \frac{N-1}2 - \nu} .
$$
From here it follows the existence of a fixed point ${\tt q}  = O( \sigma^{k - \frac{N-1}2 - \nu})$  for Problem \equ{fpN1}, and hence
$\hh = \hh^k+ {\tt q}$ satisfies constraints \equ{assh} and solves System \equ{equc}.
This concludes the proof of the theorem. \qed

\setcounter{equation}{0}
\section{Inverting the linearized Jacobi-Toda operator}
In this section we will prove Proposition \ref{solving}. The first part of the result holds in larger generality.
Actually the properties we will use in the matrix function
${\bf A}(y,\sigma)$ are  its symmetry,  its  smooth dependence
in its variables on $\KK\times [0,\sigma_0)$, and the fact that  for certain numbers  $\gamma_\pm>0$, we have
\be
\gamma_- |\xi|^2 \le \xi^T {\bf A}(y,\sigma)\xi \le \gamma_+ |\xi|^2 \foral \xi\in \R^{m-1},  \ (y,\sigma) \in \KK\times [0,\sigma_0)    .
\label{elip}\ee

\bigskip
Most of the work in the proof consists in finding the sequence $\sigma_\ell$ such that $0$
lies suitably away from the spectrum of $L_{\sigma_\ell}$,
when this operator is regarded as  self-adjoint in $L^2(\KK)$.
The result will be a consequence of various considerations on the asymptotic behavior of
the small eigenvalues of  $L_\sigma$ as $\sigma \to 0$.
The general scheme below has already been used in related settings, see \cite{mm1,mm2,mw,mahm,maz-pac0},
using the theory of smooth and analytic
dependence of eigenvalues of families of Fredholm operators due to T. Kato\cite{Kato}.
Our proof relies only on elementary considerations on the variational characterization of the
eigenvalues of $L_\sigma$ and Weyl's asymptotic formula.

\medskip
As in the above mentioned works, the assertion holds not only along a sequence, but actually for all values of $\sigma$ inside
 a sequence of disjoint intervals centered at the $\sigma_\ell$'s of width $O(\sigma_\ell^{\frac{N-1}2 })$. The corresponding assertion for
 $N=2$ can be made much more precise.

\medskip
Thus, we consider the  eigenvalue problem
\be
L_\sigma \phi = \la\phi  \quad\hbox{in }\KK.
\label{ev1}\ee
 For each $\sigma>0$ the eigenvalues are given by a sequence $\la_j(\sigma)$, characterized by the Courant-Fisher formulas
\be
\la_j(\sigma) = \sup_{ dim(M) = j-1 } \inf_{\phi\in M^\perp\setminus\{0\}}Q_{\sigma}(\phi,\phi)
=   \inf_{ dim(M) = j } \sup_{\phi\in M\setminus\{0\}} Q_{\sigma}(\phi,\phi),
\label{courant}
\ee
where
$$
Q_{\sigma}(\phi,\phi) = \frac{\int_\KK L_\sigma \phi\cdot \phi}{\int_\KK |\phi|^2} =   \frac{\int_\KK  \sigma |\nn \phi|^2 -  \phi^T {\bf A}(y,\sigma)  \phi}{\int_\KK |\phi|^2}.
$$
We have the validity of the following result.

\begin{lemma}\label{lemi0}
There is a number $\sigma_*>0$ such that for all $0< \sigma_1 < \sigma_2 < \sigma_*$ and all $j\ge 1$ the following inequalities hold.
\be
(\sigma_2 -\sigma_1) \frac {\gamma_-} {2\sigma_2^2} \le  \sigma_2^{-1}\la_j(\sigma_2)- \sigma_1^{-1}\la_j(\sigma_1)  \le  2(\sigma_2 -\sigma_1)\frac{\gamma_+}{\sigma_1^2}.
\label{desig}\ee
In particular, the functions $\sigma\in (0,\sigma_*)\mapsto \la_j(\sigma)$ are continuous.
\end{lemma}

\proof
Let us consider small numbers $0< \sigma_1 < \sigma_2$.
We observe that for any $\phi$ with $\int_\KK |\phi|^2 =1$ we have
\begin{align}
\sigma_1^{-1} Q_{\sigma_1}(\phi,\phi)  - \sigma_2^{-1} Q_{\sigma_2}(\phi,\phi) =  - \int_\KK \phi^T (\sigma_1^{-1} {\bf A}(y,\sigma_1)- \sigma_2^{-1} {\bf A}(y,\sigma_2))\phi
\nonumber
\\
\nonumber
\\
=(\sigma_1-\sigma_2) \int_\KK \phi^T (\sigma^{-2} {\bf A}(y,\sigma)
- \sigma^{-1} \partial_\sigma {\bf A}(y,\sigma))\phi,
\label{pon}
\end{align}
for some $\sigma\in (\sigma_1,\sigma_2)$.
From the assumption \equ{elip} on the matrix $A$ we then find that

$$
 \sigma_1^{-1}Q_{\sigma_1}(\phi,\phi) +(\sigma_2 -\sigma_1) \frac {\gamma_-} {2\sigma_2^2} \le  \sigma_2^{-1}Q_{\sigma_2}(\phi,\phi)  \le \sigma_1^{-1}Q_{\sigma_1}(\phi,\phi) + 2(\sigma_2 -\sigma_1)\frac{\gamma_+}{\sigma_1^2}.
$$
From here, and formulas \equ{courant}, inequality \equ{desig}  follows. \qed

\begin{corollary}\label{lemi}
There exists a number $\delta>0$ such that for any $\sigma_2>0$ and $j$
such that
  $$ \sigma_2 + |\la_j(\sigma_2) | < \delta,$$
   and  any $\sigma_1$ with $\frac12 \sigma_2 \le \sigma_1 <\sigma_2$, we have that
$$\la_j(\sigma_1) \ < \  \la_j(\sigma_2). $$

\end{corollary}

\proof
Let us consider small numbers $0< \sigma_1 < \sigma_2$ such that $\sigma_1 \ge \frac {\sigma_2} 2 . $
Then from \equ{desig} we find that
$$
\la_j(\sigma_1) \le \la_j(\sigma_2)  + (\sigma_2 - \sigma_1 ) \frac 1 {\sigma_2} [ {\la_j(\sigma_2) } -  \gamma  \frac {\sigma_1} {\sigma_2}],
$$
for some $\gamma>0$.
From here the desired result immediately follows. \qed

\subsection{Proof of Proposition \ref{solving}, general $N$ }\label{loto}

Let us consider the numbers $\bar\sigma_\ell := 2^{-\ell}$ for large $\ell \ge 1$. We will find a sequence of values $\sigma_\ell\in (\bar\sigma_{\ell+1},\bar\sigma_\ell)$ as in the statement of the lemma.

We define
\begin{align}
\digamma_\ell=\big\{\sigma\in(\bar\sigma_{\ell+1},\bar\sigma_\ell):\ \mbox{ker}L_{\sigma}\neq \{0\} \big\}.
\end{align}
If $\sigma\in \digamma_\ell$ then for some $j$ we have that $\la_j(\sigma) = 0$.
It follows that  $\la_j(\bar\sigma_{l+1}) <0$.  Indeed, let us assume the opposite. Then, given $\delta>0$, the continuity of $\la_j$ implies the existence of $\ttt\sigma$ with  $\frac 12 \sigma \le  \ttt\sigma < \sigma $ and
$0\le \la_j(\ttt \sigma) <\delta $. If $\delta$ is chosen as in Corollary \ref{lemi}, and $\ell$ is so large that  $2^{-\ell}< \delta$, we  obtain a contradiction.

\medskip
As a conclusion, we find that for all large $\ell$
\be
{\hbox{card}\,( \digamma_\ell)} \le N(\bar \sigma_{\ell+1}),
\label{cads}\ee
where $N(\sigma)$ denotes the number of negative eigenvalues of problem \equ{ev1}. We  estimate next this number for small $\sigma$.
Let us consider   $a_+>0$ such that
$$
 \xi^T {\bf A}(y,\sigma)\xi \le a_+ |\xi|^2 \foral \xi\in \R^{m-1},  \ (y,\sigma) \in \KK\times [0,\sigma_0)    ,
$$
and the operator
\begin{align}
L^+_\sigma=-\Delta_\KK -\frac{a^+}{\sigma}.
\label{k 6}
\end{align}
Let $\la^+_j(\sigma)$ denote its eigenvalues.  From the Courant-Fisher characterization we see that
$
\la^+_j(\sigma)\leq \la_j(\sigma).
$
Hence $N(\sigma) \le N_+(\sigma)$, where the latter quantity designates the number of negative eigenvalues of $L^+_\sigma$.

 \medskip
 Let us denote
 by ${\mu}_j$ the eigenvalues of $-\Delta_\KK$. Then  Weyl's asymptotic formula for eigenvalues of the Laplace-Beltrami operator, see for instance \cite{chavel,LiYau,minakPleijel}, asserts that for a certain constant $C_{\KK}>0$,
\begin{align}
\mu_j =  C_{\KK}{j}^{\frac{2}{N-1}} \ +\  o( j^{\frac{2}{N-1}})\quad\hbox{as } j\to +\infty.
\label{k 9}
\end{align}
 Using the fact that $\la^+_j(\sigma)=\mu_j-\frac{a^+}{\sigma}$  and \equ{k 9} we then find that
\begin{align}
N_+(\sigma) =  {C}{\sigma^{-\frac{N-1}{2}}}    +  o(\sigma^{-\frac{N-1}{2}})  \quad\hbox{as } \sigma\to 0 .
\label{k 10}
\end{align}
As a conclusion, using \equ{cads} we find
\be
{\hbox{card}\,( \digamma_\ell)} \le N(\bar \sigma_{\ell+1}) \le C { \bar\sigma_{\ell+1}^{-\frac{N-1}{2}}} \le C 2^{\ell \frac{N-1}{2}}.
\label{cads1}\ee
Hence there exists an interval $(a_\ell, b_\ell)\subset (\bar \sigma_{\ell+1}, \bar\sigma_\ell)$
such that $a_\ell, b_\ell\in \digamma_\ell$, $\hbox{Ker}\,(L_\sigma)=\{0\}$, $\sigma\in (a_\ell, b_\ell)$ and
\begin{align}
\label{k 17}
b_\ell-a_\ell\geq \frac{\bar \sigma_\ell-\bar \sigma_{\ell+1}}{\hbox{card}\,(\digamma_\ell)}\geq C\bar \sigma_\ell^{1+\frac{N-1}{2}}.
\end{align}

Let $$ \sigma_\ell:=\frac{1}{2}(b_\ell+a_\ell).$$ We will analyze the spectrum of $L_{\sigma_\ell}$. If some $c>0$,  and all $j$ we have
\begin{align}
|\lambda_j(\sigma_\ell)|\geq  c\bar\sigma_\ell^{\frac{N-1}2 },\label{8point12}
\end{align}
then we have  the validity of the existence assertion and estimate
(\ref{k 15}). Assume the opposite, namely that  for some $j$
we have $|\lambda_j(\sigma_\ell)|\leq  \delta \sigma_\ell^{\frac{N-1}2 }$, with
$\delta$ arbitrarily small.
 Let us assume first that
$
0< \la_j(\sigma_\ell) <  \delta \sigma_\ell^{\frac{N-1}2 }.
$
Then we have from Lemma \ref{lemi0},
$$
\la_j(a_\ell ) \le  \la_j(\sigma_\ell) - (\sigma_\ell - a_\ell ) \frac 1 {\sigma_{\ell}} [ {\la_j(\sigma_\ell ) } +  \gamma  \frac {a_\ell} { 2\sigma_\ell}].
$$
Hence, (\ref{k 17}) and (\ref{8point12}) imply that
$$
\la_j(a_\ell ) \le  \delta \sigma_\ell^{\frac{N-1}2 }
-C\frac{{\bar\sigma}_{\ell}}{2\sigma_\ell}\,{\bar\sigma}_\ell^{\frac{N-1}2 }
\Big[\lambda_j(\sigma_\ell)+\frac{\ga_{-}a_\ell}{2\sigma_\ell}\Big]<0,
$$
if $\delta$ was chosen a priori sufficiently small.
It follows that $\la_j(\sigma)$ must vanish at some $\sigma\in (a_\ell, \sigma_\ell)$,
and we have thus reached a contradiction  with the choice of the interval  $(a_\ell, b_\ell)$.

\medskip
The case $-\delta \sigma_\ell^{\frac{N-1}2 }< \la_j(\sigma_\ell) < 0$ is handled similarly. In that case we get $\la_j(b_\ell) >0$.
 The proof of existence and estimate  (\ref{k 15}) is thus complete.

\medskip
Let us consider now a number $p>N-1$. Now we want to estimate the inverse of $L_{\sigma_\ell}$ in Sobolev norms.
The equation satisfied by $\psi = L_{\sigma_\ell}^{-1} g$ has the form
$$
\Delta_{\KK} \psi  = O(\sigma^{-1})[ \psi +{{g}} ]
$$
for $\sigma =\sigma_\ell$.
Then from elliptic estimates we get
\be
 \|\psi\|_{W^{2,q}(\KK)} \le C\sigma^{-1} [ \|\psi\|_{L^q(\KK)} + \|g\|_{L^q(\KK)} ]
\label{er}\ee
Using this for $q=2$ and estimate  (\ref{k 15})  we obtain

$$ \|\psi\|_{W^{2,2}(\KK)} \le C\sigma^{-1} [ \|\psi\|_{L^2(\KK)} + \|g\|_{L^2(\KK)} ] \le  C\sigma^{- \frac{N}2-1} \|g\|_{L^p(\KK)}.
$$
From Sobolev's embedding we then find

$$ \|\psi\|_{L^q(\KK)} \le   C\sigma^{- \frac{N-1}2-1} \|g\|_{L^p(\KK)}.
$$
for any $1< q \le  \frac {2(N-1)}{N-5}$ if $N> 5$, and any $q>1$ if $N\le 5$. If $q=p$ is admissible in this range, the estimate follows from \equ{er}. If not,
we apply it for  $q=   \frac {2(N-1)}{N-5}$, and then Sobolev's embedding yields

$$ \|\psi\|_{L^s(\KK)} \le   C\sigma^{- \frac{N-1}2-2 } \|g\|_{L^p(\KK)}.
$$
for any $1< s \le  \frac {2(N-1)}{N-6}$ if $N> 6$, and any $s>1$ if $N\le 6$. Iterating this argument, we obtain the desired estimate after a finite number of steps.  The proof of the first part of the proposition is  concluded. \qed

\subsection{The case $N=2$. Conclusion of the proof }
We consider now the problem of solving system \equ{bvp} when $N=2$.
We consider first the problem of solving

\be
- \sigma \Delta_\KK \psi - {\bf A}(y,0)\psi\ =\  g   \quad\hbox{in }\KK.
\label{bvp1-1}
\ee

A main observation is the following:  the linear system \equ{bvp1-1} can be decoupled:
If  $\Lambda_1, \ldots, \Lambda_{m-1}$ denote the eigenvalues of the matrix
$$
{\bf Q}:= {\bf C}^{\frac 12}\,\left [\begin{matrix}   a_1 &  0& \cdots & 0
\\
                                                                                0 &   a_2& \cdots & 0
\\
                                                                                 \vdots &                   & \ddots &  \vdots  \\
                                                                                0   &  0& \cdots & a_{m-1}
                                                                                \end{matrix} \right ]\, {\bf C}^{\frac 12}, $$
                                                                                which coincide with those of
$$
 \left[\begin{matrix} 2a_1 & -a_2&  0 & &\cdots  &  0 \\
-a_1 & 2a_2& -a_3 & &\cdots   &  0  \\
                       0    &   -a_2 &  2a_3 & &\cdots & 0\\
                      \vdots   &  &\ddots & \ddots & \ddots   & \vdots   \\
                       0 &   & \cdots   & -a_{m-3} & 2a_{m-2} &  -a_{m-1} \\
                       0 &    &\cdots & & -a_{m-2}&  2a_{m-1} \\
                          \end{matrix}\right ],
                                                                                $$
then system \equ{bvp1-1} expressed in  coordinates associated to eigenfunctions of ${\bf Q}$ decouples into
$m-1$ equations of the form
\be
- \sigma \Delta_\KK \psi_j  - \frac \beta{\sqrt{2}} \Lambda_j\,K(y)\, \psi_j\ =\  g_j
 \quad
\hbox{in }\KK, \quad j=1,\ldots, m-1.
\label{bvp1}
\ee

When $N=2$ this problem reduces to an ODE. $\KK$ is then a geodesic of $\MM$ and $K(y)$ will  simply be Gauss curvature measured along $\KK$.
Using $y$ as arclength coordinate, and dropping the index $j$, Equations \equ{bvp1} take the generic form

\begin{align}
\begin{aligned}
- \sigma \psi'' -   \mu\, K(y)\,\psi =  g \quad \mbox{in } (0,\ell),
\\
\psi(0)\,=\,\psi(\ell ),\quad\psi'(0)\,=\,\psi'(\ell ),
\end{aligned}
\label{eigenvalueproblem2}
\end{align}
where $\mu$ is given and fixed, and
 $\ell$ is the total length of $\KK$.

For this problem to be uniquely solvable, we need that $\mu\sigma^{-1}$ differs  from the eigenvalues $\la= \la_j$ of
 the  problem

 \begin{align}
\begin{aligned}
-  \varphi'' =   \la \, K(y)\,\varphi  \quad \mbox{in } (0,\ell),
\\
\varphi(0)\,=\, \varphi(\ell ),\quad \varphi'(0)\,=\,\varphi'(\ell ).
\end{aligned}
\label{eigenvalueproblem}
\end{align}
More precisely, in such a case we have that the solution of \equ{eigenvalueproblem2} satisfies
\be
\|\psi\|_{L^2(\KK)} \le  \frac {\sigma^{-1}}{\min_j |\la_j - \sigma^{-1}\mu|} \|g\|_{L^2(\KK)}.
\label{pico}\ee

Now, we restate Problem \equ{eigenvalueproblem}
using the following Liouville transformation:
\bas
 \ell_0=\int^{\ell }_0\sqrt{K(y)}\,\mathrm{d}\,y,\qquad
t=\frac{\pi}{\ell_0}\int^{y}_0\sqrt{K(\theta)}\,\mathrm{d}\,\theta,\, t\in [0,\pi),
\\
\Psi(y)=K( y)^{-\frac{1}{4}},\qquad e(t)=\varphi(y)/\,\Psi(y),\qquad q(t)=\frac{\ell_0^2\,\Psi''(y)}{\pi^2 \,\Psi^2(y)\,K(y)}.
\eas
Equation \equ{eigenvalueproblem} then becomes
\bas
-e''-q(t)\,e=\frac{\ell_0^2}{\pi^2}\,\la\,e\quad \mbox{ in } (0,\pi),\
e(0)\,=\,e(\pi),\ e'(0)=e'(\pi).
\eas

A result in \cite{LS} shows that, as $j\rightarrow \infty $ we have
$$
\la_j = \frac{4\pi^2j^2}{\ell_0^2}  \,+\,O(j^{-2}).
$$
Hence, if for some $c>0$ we have that

$$
\Big|\sigma^{-1} \mu - \frac{4\pi^2j^2}{\ell_0^2}\Big| \,>\, c\sigma^{- \frac 12} \foral j\ge 1,
$$

\noindent and $\sigma$ is sufficiently small, then the problem will be solvable, and thanks to \equ{pico}, we will have the estimate

\be
\|\psi\|_{L^2(\KK)} \le   C\sigma^{-\frac 12} \|g\|_{L^2(\KK)},
\label{estt}\ee

\noindent
for the unique solution of   Problem \equ{eigenvalueproblem2}. It follows that, under these conditions System \equ{bvp1} is uniquely solvable
and its solution $\psi = - (\sigma \Delta_\KK \psi + {\bf A}(y,0))^{-1}g $ satisfies estimate \equ{estt}.

\medskip
Now, for $\sigma$ as above, we can write system \equ{bvp} in the fixed point form in $L^2(\KK)$,

\be
\psi  + T(\psi) =   - (\sigma \Delta_\KK \psi + {\bf A}(y,0))^{-1}g  ,\quad  \psi\in L^2(\KK),
\label{pata}\ee

\noindent where

$$
T(\psi):=   (\sigma \Delta_\KK \psi + {\bf A}(y,0))^{-1} [ ({\bf A}(y,\sigma) -{\bf A}(y,0))\psi ]\, .
$$

\noindent We observe that, as an operator in $L^2(\KK)$,  $\|T\| = O(\sigma^{\frac 12})$. Thus, for small $\sigma$, Problem \equ{pata} is
uniquely solvable, and satisfies \equ{estt}.  Finally, for the $L^p$ case, we argue with the same bootstrap procedure of \S \ref{loto}.

\medskip
The proof of the proposition is complete. \qed
\bigskip

\noindent {\bf Acknowledgments.}
The first and second authors have been supported by grants Anillo ACT 125 Center for Analysis of Partial Differential Equations (CAPDE), FONDECYT 1070389 and 1090103,
and by Fondo Basal CMM-Chile.
The third author is supported by an Earmarked Grant (GRF) from RGC of Hong Kong
and a ``Focused Research Scheme'' from CUHK.
The fourth author is supported by grants 10571121 and 10901108 from NSFC.
He thanks the department of Mathematics of the Chinese University of Hong Kong for its kind hospitality.
We thank  Frank Pacard for valuable remarks.

\end{document}